\documentclass[11pt]{article}

\usepackage{amssymb,amsmath,amsfonts,amssymb,numbysec}
\usepackage{graphics,graphicx,color,hhline}
\usepackage{bbm}
\usepackage{euscript} 

\def\@abssec#1{\vspace{.05in}\footnotesize \parindent .2in 
{\bf #1. }\ignorespaces} 

\graphicspath{{/EPSF/}{Figures/}}   

\newtheorem{theorem}{Theorem}[section]
\newtheorem{lemma}[theorem]{Lemma}
\newtheorem{proposition}[theorem]{Proposition}

\def \Rm {\mathbb R}

\newcommand{\eps}{\varepsilon}
\newcommand{\E}{\mathbb E}

\newcommand{\pdr}[2]{\dfrac{\partial{#1}}{\partial{#2}}}

\newcommand{\bn}{\mathbf n}
  
\newcommand{\bu}{\mathbf u} \newcommand{\bv}{\mathbf v}

\newcommand{\calQ}{\mathcal Q}

\newcommand{\calF}{\mathcal F}
\newcommand{\calP}{\mathcal P}
\newcommand{\calO}{\mathcal O} 

\newcommand{\Tr}{\textrm{Tr}}

\newcommand{\calV}{\mathcal V}

\newcommand{\calS}{\mathcal S}

\newcommand{\frakH}{\mathfrak H}

\newcommand{\bzero}{\mathbf 0}
\newcommand{\eC}{\EuScript{C}}
\def\fref#1{{\rm (\ref{#1})}}

\newcommand{\cout}[1]{}

\def\un{{\mathbbmss{1}}}

\newcommand{\pr}{\parallel}

\newcommand{\be}{\begin{equation}}
\newcommand{\ee}{\end{equation}}
\newcommand{\bea}{\begin{eqnarray}}
\newcommand{\eea}{\end{eqnarray}}
\newcommand{\bee}{\begin{eqnarray*}}
\newcommand{\eee}{\end{eqnarray*}}
\newcommand{\bal}{\begin{align*}}
\newcommand{\eal}{\end{align*}}
\numberbysection

\begin{document}
{\title{On the stability of some imaging functionals}}
\author{Guillaume Bal  \thanks{Department of Applied Physics and
     Applied Mathematics, Columbia University, New York NY, 10027 ;
     gb2030@columbia.edu}
\and Olivier Pinaud \thanks{Department of Mathematics,
     Colorado State University, Fort Collins CO, 80523 ;
pinaud@math.colostate.edu} \and
Lenya Ryzhik \thanks{Department of Mathematics,
     Stanford University, Stanford CA 94305 ; ryzhik@math.stanford.edu}}

\maketitle

\begin{abstract}
  This work is devoted to the stability/resolution analysis of several
  imaging functionals in complex environments. We consider both linear
  functionals in the wavefield as well as quadratic functionals based
  on wavefield correlations. Using simplified measurement settings
  and reduced functionals that retain the main features of functionals
  used in practice, we obtain optimal asymptotic estimates of the signal-to-noise
  ratios depending on the main physical parameters of the problem. We
  consider random media with possibly long-range dependence and with a 
  correlation length that is less than or equal to the central wavelength
  of the source we aim to reconstruct. This corresponds to the wave propagation regimes of  
  radiative transfer or homogenization.
\end{abstract}

\section{Introduction}  
Imaging in complex media has a long history with many applications
such as non-destructive testing, seismology, or underwater acoustics
\cite{Claerbout, Shull, bag-math}. Standard methodologies consist in
emitting a pulse in the heterogeneous medium and performing
measurements of the wavefield or other relevant quantities at an array of
detectors. Depending on the experimental setting, this may give
access, for instance, to the backscattered wavefield (and its
spatio-temporal correlations), or to the wave energy. The imaging
procedure then typically amounts to back-propagating the measured data
appropriately. When scattering of the wave by the medium is not too
strong, measurements are usually migrated in a homogeneous medium
neglecting the heterogeneities. This is the principle of the Kirchhoff
migration and similar techniques, and is referred to as coherent
imaging~\cite{bleistein}. When scattering is stronger, the interaction
between the wave and the medium cannot be neglected and a different
model for the inversion is needed. We will only consider weak
heterogeneities in this work, so that a homogeneous model is enough to
obtain accurate reconstructions. We refer to \cite{BR-AM,BP-M3AS} for
a consideration of stronger fluctuations and transport-based imaging.

We are interested here in comparing two classes of imaging methods in
terms of stability and resolution: one based on the wavefield
measurements (and thus linear in the wavefields), such as Kirchhoff migration; and the
other one based on correlations of the wavefield (and therefore
quadratic in the wavefields), such as coherent interferometry imaging \cite{CINT}. 
Here, ``stability'' refers to the stability of the reconstructions with respect to
changes in the medium or to a measurement noise. 

Stability and
resolution of coherent imaging functionals in random waveguides were
addressed in \cite{Borcea-wg}; for functionals based on topological
derivatives, see also e.g. \cite{habib-topo}. 
We refer to
\cite{FGPS-07,Garnier-P, GS-velo, CINT-SNR} and references therein
for more details on wave propagation in random media and imaging.
Our analysis is done in the framework of three-dimensional acoustic
wave propagation in a random medium with correlation length $\ell_c$
and fluctuations amplitude of order~$\sigma_0$. The quantities
$\ell_c$ and $\sigma_0$ are fixed parameters of the experiments. The
random medium is allowed to exhibit either short-range or long-range
dependence. Our goal is to image a source centered at the origin from
measurements performed over a detector array $D$. Denoting by $\lambda$ the central
wavelength of the source, we will assume that it is larger than or
equal to the correlation length, and that many wavelengths separate
the source from the detector, so that we are in a high
frequency regime. The case $\lambda=\ell_c$ is referred to as the weak
coupling regime in the literature (for an amplitude $\sigma_0\ll 1$),
while the case $\lambda \gg \ell_c$ is the stochastic homogenization
regime. The ratio $\ell_c/\lambda$ is a critical parameter as it
controls, along with $\sigma_0$, the strength of the interaction of
the sound waves with the heterogeneities. The regime $\lambda \ll \ell_c$ is
addressed in~\cite{CINT-Habib} and is known as the random
geometrical optics regime. The stability/resolution analysis is
somewhat direct in that case, since a simple expression for the
heterogeneous Green's function is available by means of  
random travel times. This is not the case in our regimes of interest,
where the interaction between the waves and the underlying medium is 
more difficult to describe mathematically.

We will work with a simplified configuration and will define
reduced imaging functionals that retain the main features of the
functionals used in practice while offering more tractable computations. In
such a setting, we obtain optimal estimates for the stability of the
functional in terms of the most relevant physical parameters. We
furthermore quantify the signal-to-noise ratio (SNR). 
The statistical
stability is evaluated by computing the variance of the functionals at
the source location. For the wave-based functional (WB), this involves
the use of stationary phase techniques for oscillatory integrals,
while computations related to the correlation-based functional (CB)
involve averaging the scintillation function of Wigner transforms
\cite{GMMP,LP} against singular test functions. It is now
well-established that correlation-based functionals enjoy a better
stability than wave-based functionals at the price of a lower
resolution \cite{CINT}. This is due to self-averaging effects of
Wigner transforms that we want to quantify here.

The paper is organized as follows: in Section \ref{mainresults}, we
present the setting and our main results. We define the wave-based and
correlation-based functionals and describe our models for the
measurements in Section \ref{meas_imag}. Section \ref{secCB} is
devoted to the analysis of the resolution and stability of the
correlation-based functional and Section \ref{secWB} to that of the
wave-based functional. The proofs of the main results are presented in
Section \ref{secproofs} and a conclusion is offered in Section
\ref{conc}.
%
%

\section{Setting and main results} \label{mainresults}

\subsection{Wave propagation and measurement setting} \label{gene}

The propagation of three-dimensional acoustic waves is described by
the  scalar
wave equation for the pressure $p$:
 $$ \frac{\partial^2 p}{\partial t^2}=\kappa^{-1}(x) 
   \nabla \cdot [\rho^{-1}(x)  \nabla p], \qquad x \in \Rm^3, \;t>0,
$$
supplemented with initial conditions $p(t=0,x)$ and $p_t(t=0,x)$. We
suppose for simplicity that the density is constant and equal to one;
$\rho=1$. 
We also assume that the compressibility is random and takes the form
$$\kappa(x)=\kappa_0 \left(1+\sigma_0 V \left(\frac{x}{\ell_c} \right)\right), \qquad \kappa_0=1,$$
where $\sigma_0$ measures the amplitude of the random fluctuations and
$\ell_c$ their correlation length. We suppose that $V$ is bounded and
$\sigma_0$ is small enough so that $\kappa$ remains positive. The
sound speed $c=(\kappa\rho)^{-\frac12}$ thus satisfies $c=1+\calO(\sigma_0) \simeq 1$
and the average sound speed is $c_0=1$. The
random field $V$ is a mean-zero stationary process with correlation
function
$$R(x)=\E \{V(x+y) V(y)\},$$
where $\E$ denotes the ensemble average over realizations of the
random medium. We assume~$R$ to be isotropic so that $R(x)=R(|x|)$. We
will consider two types of correlation functions: (i) integrable $R$,
which models random media with short-range correlations; and (ii)
non-integrable $R$ corresponding to media with long-range
correlations. The latter media are of interest for instance when waves
propagate through a turbulent atmosphere or the earth upper crust
\cite{dolan,sidi}. Such properties can be translated into the power
spectrum $\hat{R}$, the Fourier transform of $R$, by supposing that
$\hat{R}$ has the form
$$
\hat{R}(k)=\frac{S(k)}{|k|^\delta}=\int_{\Rm^3} e^{-i x \cdot k} R(x)dx, \qquad 0 \leq \delta <2,
$$
where $S$ is a smooth function with fast decay at infinity. A simple dimensional analysis shows that $R(x)$ behaves likes $|x|^{\delta-3}$, $\delta<3$, as $|x|\to\infty$, which is not integrable for large $|x|$. The case $\delta=0$ corresponds to integrable $R$ since the function $S$ is regular. We assume that $\delta<2$ so that the transport mean free path (see section \ref{models}) is well-defined. Propagation in media such that $2<\delta<3$ is still possible but requires an elaborate theory of transport equations with highly singular scattering operators that are beyond the scope of this paper.

\begin{figure}[here]
   \begin{center}
\includegraphics[height=6cm, width=11.2cm]{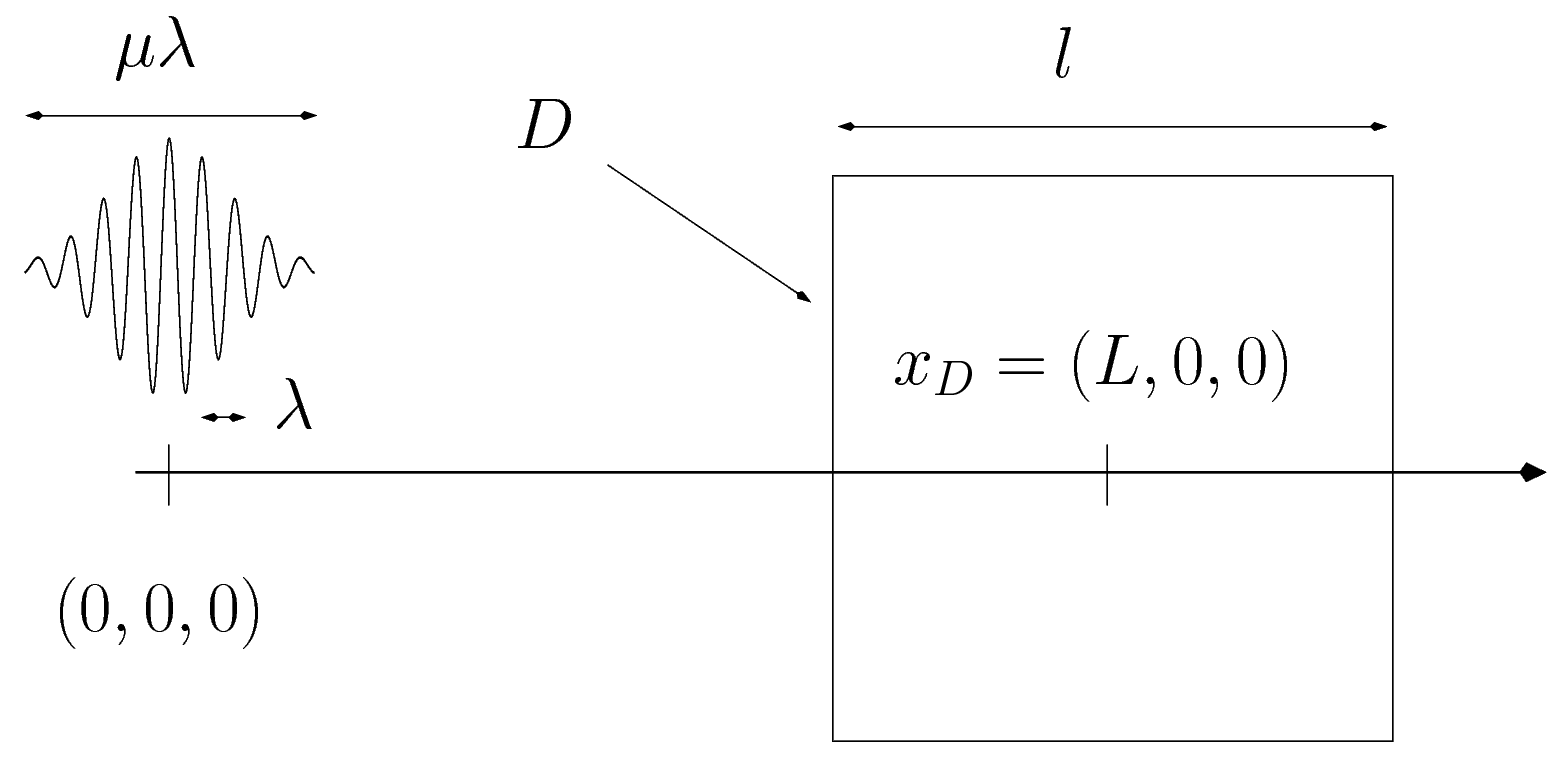}
\end{center}
\caption{Geometry}
\label{fig1}
\end{figure}

Our setting of measurements is the following: we assume that measurements of the wavefield are performed on a three dimensional detector $D$ and at a fixed time $T$. We assume that the detector $D$ is a cube centered at $x_D=(L,0,0)$ of side $l<L$; see figure \ref{fig1}.  Our goal is then to image an initial condition 
$$p(t=0,x)=p_0(x), \qquad \partial_t p(t=0,x)=0,$$
localized around $x_0 =(0,0,0)$.  

Practical settings of measurements usually involve recording in time of the wavefield on a surface detector. The two configurations share similar three-dimensional information: spatial 3D for our setting, and spatial 2D + time 1D for the standard configuration. Using wave propagation in a homogeneous medium, it is also relatively straightforward to pass from one measurement setting to the other. Our choice of a 3D detector was made because it offers slightly more tractable computations for the stability of the imaging functionals while qualitatively preserving the same structure of data. 

As waves approximately propagate in a homogeneous medium with constant speed $c_0$, it takes an average time $t_D=(L-l/2)/c_0$ for the wave to reach the detector. We will therefore suppose that $T>t_D$, and even make the assumption that $T= L/c_0$ so that the wave packet has reached the center of the detector. Note that in such a measurement configuration, the range of the source is then  known since the initial time of emission is given. We therefore focus only on the cross-range resolution.

We consider an isotropic initial condition obtained by Fourier transform of a frequency profile $g$ (a smooth function that decays rapidly at infinity):
$$
p_0(x)= \lambda^3 \int_{\Rm^3} g( B^{-1} (|k|-|k_0|/\lambda)) e^{-i k \cdot x}  dk,
$$
where $|k_0|$ is a non-dimensional parameter, $\lambda$ is the central wavelength, and $B c_0^{-1}$ the bandwidth. 

Rescaling variables as $x \to x L$, $t \to t L/c_0 $, introducing the parameters $\eta=\ell_c/\lambda$, $\eps=\lambda/L$, and still denoting by $p$ the corresponding rescaled wavefield, the dimensionless wave equation now reads
\be \label{waveq}
\left(1+\sigma_0 V \left(\frac{x}{\eta \eps} \right) \right) \frac{\partial^2 p}{\partial t^2}=\Delta  p, \quad p(t=0,x)=p_{0}^\eps(x), \quad \frac{\partial
  p}{\partial t}(t=0,x)=0.
\ee
We quantify the bandwidth of the source in terms of the central frequency by setting $B=(L \eps \mu)^{-1}$, where $\mu$ is a given non-dimensional parameter such that  $\mu \gg 1$. We suppose that $\mu \ll 1/\sqrt{\eps}$ so that the initial condition can be considered as broadband. More details about the latter condition will be given later. The rescaled initial condition then has the form
\be \label{CI1}
p_{0}^\eps(x)=\mu \eps^{3} \int_{\Rm^3} g\left( \eps \mu\left(|k|-\frac{|k_0|}{\eps}\right)\right) e^{-i k \cdot x}  dk.
\ee
The normalization is chosen so that $p_0^\eps(0)$ is of order $O(1)$. If $\hat g$ denotes the Fourier transform of $g$, we can write
$$
p_{0}^\eps(x) = p_{0}^\eps(|x|)\simeq |k_0|^2 \int_{S^2} \hat g\left(\frac{\hat k \cdot x}{ \eps \mu}\right)e^{-i \frac{|k_0|}{\eps} \hat k \cdot x} d \hat k.
$$
Above, the symbol $\simeq$  means equality up to negligible terms and $S^2$ is the unit sphere of $\Rm^3$. The initial condition is therefore essentially a function with support $\eps \mu$ oscillating isotropically at a frequency $|k_0|/\eps$. 
As explained in the introduction, we assume  that $\ell_c \leq \lambda$, which implies $\eta \leq 1$. The case $\eta \sim 1$ leads to the radiative transfer regime in the limit $\eps \to 0$ \cite{RPK-WM}. The case $\eta\ll1$ corresponds to the homogenization regime \cite{Bal-homog} and a  propagation in an effective medium (which is here the homogeneous medium since $\sigma_0$ is small). The opposite case $\eta \gg 1$ gives rise to the random geometric regime addressed in \cite{CINT-Habib}. 

For the asymptotic analysis, we recast the scalar wave equation (\ref{waveq}) as a first-order hyperbolic system on the wavefield $\bu=(\bv, p)$, where $\bv$ is the velocity:
$$
\rho \frac{\partial \bv}{ \partial t }+\nabla p=0,
 \qquad \kappa \frac{\partial p}{ \partial t }+\nabla \cdot \bv=0,
$$
augmented with initial conditions $p(t=0,x)=p_0^\eps\left(x \right)$ and
$\bv(t=0,x)=0$. The latter system is rewritten as
\begin{equation}
  \label{eq:hypsyst}
  \left( I+\sigma_0 \calV \left(\frac{x}{\eps \eta}\right)\right)\pdr{\bu}t + D^j \pdr{\bu}{x_j} =0,
\end{equation}
where $\calV=\textrm{diag}(\bzero, V) \in \Rm^{4 \times 4}$, $D^j_{mn}=\delta_{m4}\delta_{nj}+
\delta_{n4}\delta_{mj}$ is a $4\times4$ symmetric matrix for $1\leq
j\leq3$.  Above and in the sequel, we use the Einstein convention of summation over repeated
indices.

\subsection{Results} \label{results}
Our main results concern the signal-to-noise ratio at the point $x=0$ defined by
$$
SNR=\frac{\E \{I\}(0) }{\sqrt{\textrm{Var} \{ I\}(0)}}, \qquad I=I^C \textrm{  or  } I^W,
$$
where $I^C$ stands for the CB functional and $I^W$ for the WB functional, which are defined further in section \ref{meas_imag}. We will also quantify the support of $\textrm{Var} \{ I\}(z)$. We will obtain optimal estimates for the CB and WB functionals in terms of the main parameters that define them as well as $\eps$, $\eta$, $\delta$ and $\mu$. Our measurements of the pressure are assumed to have the form

$$p(t,x)=\E \{p\}(t,x) +\delta p^\eps (t,x)+\sigma_n n_p^\eps(t,x),$$
where $\delta p^\eps (t,x)$ models the statistical instabilities in the Born approximation and $\sigma_n n_p^\eps(t,x)$ is an additive mean zero noise with amplitude $\sigma_n$. Note that by construction, coherent-based functionals will perform well only in the presence of  relatively weak heterogeneities. Hence, modeling the statistical instabilities at first order is both practically relevant and mathematically feasible. Taking into account second order interactions is considerably more difficult mathematically; see \cite{BP-CPDE2} for a stability analysis of Wigner transforms (and not of the more complicated imaging functionals) in the paraxial approximation. The term $\sigma_n n_p^\eps$ models an additive noise at the detector and also takes into account in a very crude way the higher order interactions between the wave and the medium that are not included in $\delta p^\eps$. We suppose that $n_p^\eps$ oscillates at a frequency $\eps^{-1}$ (so that a simple frequency analysis cannot separate the real signal from the noise) and that it is independent of the random medium. The variance can then be decomposed as $\textrm{Var} \{ I^C\}=V^C+V_n^C$, where $V^C$ is the variance corresponding to the single scattering term and $V_n^C$ the noise contribution. We also write $\textrm{Var} \{ I^W\}=V^W+V_n^W$. Note that the measurement noise can actually be much larger than the average $\E \{p\}$. Indeed, a simple analysis of $\E \{p\}(t,x)$ shows that it is of order $\eps$ when $|x|=1$ (omitting the absorption). It is the refocusing properties of the functionals that lead to a reconstructed source of order one from measurements of order $\eps$. 

We formally rescale the single scattering instabilities by $e^{-c_0 \Sigma  t/2}$, where $\Sigma^{-1}$ is the transport mean free path defined in section \ref{models}, so that the first-order interaction between the average field $\E \{p\}$ (proportional to $e^{-c_0 \Sigma t/2}$) and the medium has a comparable amplitude to $\E \{p\}$. The fact that the single scattering instabilities are exponentially decreasing as the ballistic part can be proved in simplified regimes of propagation, where a closed-form equation for their variance can be obtained; see e.g. \cite{B-Ito-04, BP-CPDE}. 

We need to introduce another important parameter,  which is the typical length over which correlations are calculated in the CB functional. In dimensionless units, we define it in terms of the wavelength by $N_C \eps$. We assume that the detector cannot perform subwavelength measurements, which implies $N_C \geq 1$. We suppose for simplicity that correlations are calculated isotropically. Accounting for anisotropic correlations would add additional parameters and technicalities, but presents no conceptual difficulties. Since the resolution and the stability of the CB functional are mostly influenced by $N_C \eps$ and not by the size of the detector as $N_C \eps \leq l/L$ by construction (see  section \ref{meancint} below), we will systematically suppose that the detector is sufficiently large compared to the wavelength so that its effects on the stability can be neglected in a first approximation. The parameter $N_C$ is crucial in that it controls the resolution of the CB functional (shown in section \ref{meancint} to be $L/N_C$) and its stability. Small values of $N_C$ yield a good stability for a poor resolution, while large values lead to less stability with a resolution comparable to that of the WB functional. 
\begin{table}[ht]
\caption{Notation}
\begin{tabular}{c|lc}
$L$ & Source-detector distance; see Fig. \ref{fig1}\\
$l$ & Size of the array; see Fig. \ref{fig1}\\
$\sigma_0$ & Amplitude of the random fluctuations \\
$\ell_c$ & Correlation length of the random fluctuations\\
$\delta$ & The correlation function decreases as $|x|^{\delta-3}$, $0 \leq \delta <2$ \\
$\lambda$ & Central wavelength of probing signal \\
$\eps$, $\eta$ &  $\eps=\frac\lambda L$ with $L$ distance of propagation; $\eta=\frac{\ell_c}{\lambda}$ \\
$B$ & $Bc_0^{-1}$ is bandwidth with $c_0$ background sound speed \\
$\mu$ & $\mu=\frac1{L\eps B}$ with width of probing signal equal to $\mu\lambda$\\
$\sigma_n$ & Amplitude of detector noise $\sigma_n n^\eps_p(t,x)$\\
$I^C$ and $I^W$ & Correlation-based (CB) and wave-based (WB) imaging functionals, respectively \\
$V^C$ and $V^W$ & Variance of the  above imaging functionals\\
$N_C$ & $\eps N_C$ is length over which correlations are considered in $I^C$ \\
$\frac{\lambda L}l$ & Cross-range resolution of wave-based functional (Rayleigh formula)\\
$\frac{L}{N_C}$ & Cross-range resolution of correlation-based functional \\
$\Sigma$ & Inverse of transport mean free path; see \eqref{RTE2}.
\end{tabular}
\end{table}

Our main results are the following:
\begin{theorem} \label{theo} Denote by $V^C(z)$ (resp. $V_n^C$) and $V^{W}(z)$ (resp. $V_n^W$) the variances of the correlation-based and wave-based functionals defined in section \ref{meas_imag} for the single scattering contribution (resp. the noise contribution). Let $\Sigma^{-1}$ be the transport mean free path defined in \eqref{RTE2} in section \ref{models} below. Then we have:
\begin{align*}
&V^C(0) \sim  e^{-2\Sigma L}\sigma_0^2 \; \mu^{-2}\left(\frac{\ell_c}{\lambda}\right)^{3-\delta} \left(\frac{\lambda}{L}\right)^{4+(1-\delta)\wedge 0} \bigg( \left(\frac{L}{\lambda}\right)^4\wedge (N_C\mu)^4 \bigg)\\
&V^C(0) \sim  e^{-2\Sigma L}\sigma_0^2 \;  \left(\frac{\lambda}{L}\right)^{4} \left[\frac{1}{\mu^2}\left(\frac{\ell_c}{\lambda}\right)^{3-\delta}\left(\frac{\lambda}{L}\right)^{(1-\delta)\wedge 0}\bigg( \left(\frac{L}{\lambda}\right)^4\wedge (N_C\mu)^4 \bigg) \right] \vee \left[ \left(\frac{\ell_c}{\lambda}\right)^{(3-\delta)\wedge 2} N_C^4\right]\\
&V^C_n(0) \sim e^{-\Sigma L} \sigma_n^2 \left(\frac{\lambda}{L}\right)^{4} N_C^4 \\
&V^W(0) \sim  e^{-\Sigma L}\sigma_0^2, \qquad V^W_n(0) \sim \sigma_n^2.
\end{align*}
Moreover, $V^C$, $V^C_n$ and $V^W$ are mostly supported on $|z|\leq \mu \eps$. Above, the notation $\sim$ means equality up to multiplicative constants independent of the main parameters, and $a \wedge b=\min(a,b)$, $a \vee b=\max(a,b)$ . Denote now by $SNR^C$ (resp. $SNR_n^C$) and $SNR^{W}$ (resp. $SNR_n^W$) the corresponding signal-to-noise ratios. Then:
$$
SNR^C_{tot}(0) =  SNR^C(0)\wedge  SNR_n^C(0) \qquad SNR_{tot}^W(0) =  SNR^W(0)\wedge  SNR_n^W(0),
$$
where
\begin{align*}
&SNR^C(0) \sim  \frac{1}{\sigma_0} \left[\mu \left(\frac{\lambda}{\ell_c}\right)^{(3-\delta)/2}\left(\frac{L}{\lambda}\right)^{(((1-\delta)/2)\wedge 0)} \right] \wedge
\left[\left(\frac{\lambda}{\ell_c}\right)^{((3-\delta)/2)\wedge 1} \left(\frac{L}{N_C \lambda}\right)^2 \wedge \mu^2  \right]\\ 
& SNR^C_n(0) \sim \frac{e^{-\frac12\Sigma L}}{\sigma_n} \left(\left(\frac{L}{N_C \lambda}\right)^2 \wedge \mu^2\right) \\
&SNR^W(0) \sim  \frac{1}{\sigma_0}, \qquad SNR^W_n(0) \sim \frac{e^{-\frac12\Sigma L}}{\sigma_n}.
\end{align*}
\end{theorem}
We below comment on the results of the theorem.
\paragraph{Comparison CB-WB.} Let us consider first the case of short-range correlations $\delta=0$. We neglect absorption at this point (i.e. $e^{-\Sigma L}=1$) and will take it into account when considering the SNR of the external noise. 

Let us start with the single scattering contributions $SNR^C$ and $SNR^W$. Suppose first that $N_C$ is small enough and $\mu$ large enough so that $(L/N_C\lambda)\wedge \mu=\mu$. Then,
$$SNR^C \sim  \frac{1}{\sigma_0} \; \left(\frac{\lambda}{\ell_c} \right)\mu .$$
Hence, the SNR increases when the following quantities increase: $\lambda/ \ell_c$ (the dynamics gets closer to the homogenization regime where the homogenized solution is deterministic),  $\mu$  (smaller bandwidth, which results in a loss of range resolution) and $\sigma_0^{-1}$ (weaker fluctuations). The fact that the SNR increases as $\mu$ stems from the self-averaging effects of the quadratic functional. There are no such effects for $V^W$ and we find $SNR^C(0)=SNR^W(0) (\mu \lambda/\ell_c)$ so that the CB functional is more stable than the WB functional in this configuration of small $N_C$, whether in the radiative transfer regime ($\lambda=\ell_c$) or in the homogenization regime ($\lambda \gg \ell_c$).

It is shown in section \ref{meancint} that the resolution of the CB functional is $L/N_C$, so that the best resolution is achieved for the largest possible value $N_C$, where correlations are calculated over the largest possible domain, namely $N_C = l/\lambda$. Hence, the best resolution is $\lambda L/l$, which is the celebrated Rayleigh formula and the same cross-range resolution as the WB functional. Now, for this best resolution and when $l \mu > L$, then $(L/N_C\lambda)\wedge \mu=(L/N_C\lambda)$ so that $SNR^C(0)=SNR^W(0) (\lambda/\ell_c)$ (equality here means equality up to multiplicative constants independent of the parameters $\lambda$, $N_C$, $\ell_c$, $\sigma_0$). This means that in the absence of external noise, and in a weak disorder regime where multiple scattering can be neglected, the SNR of the CB and WB functionals differ only by the factor $(\lambda / \ell_c)$ for a similar resolution. This is a term of order one in the radiative transfer regime, and a large term in the homogenization regime. In the radiative  transfer regime, in order to significantly increase the CB SNR compared to WB, one needs to decrease $N_C$ and therefore to lower the resolution. This is the classical stability/resolution trade off as was already observed in \cite{CINT-Habib} for the random geometrical regime. Note also that the statistical errors for both functionals are essentially localized on the support (of diameter $\mu \lambda$) of the initial source.


We consider now the noise contributions $V_n^C$ and $V^W_n$. A first important observation is that $V^C_n(x)$ is mostly localized on the support of the source, while $V^W_n(x)$ is essentially a constant everywhere. This stems from the fact that in the CB functional, the noise is correlated with the average field so that the functional keeps track of the source direction, while in the WB functional the noise is correlated with itself, see sections \ref{flucmod} and \ref{secWB}. Moreover,the SNR are
$$
SNR^W_n(0) \sim \frac{e^{-\frac12\Sigma L}}{\sigma_n}, \qquad SNR^C_n(0) \sim SNR^W_n(0)\left(\left(\frac{L}{N_C \lambda}\right)\wedge \mu \right)^2.
$$
 As for the single scattering SNR, for identical resolutions (i.e. when $N_c=l/\lambda$), the SNR are comparable. In order to obtain to better SNR for CB, one needs to lower the resolution (and therefore decrease $N_C$) so that $(L/(N_C \lambda)) \wedge \mu =\mu$  and $SNR^C_n(0)=SNR^W_n(0) \mu \gg SNR^W_n(0)$. 
\paragraph{Minimal wavelength for a given SNR.} We want to address here the question of the lowest central wavelength $\lambda_m$ of the source (and therefore the best resolution) that can be achieved for a given SNR. In order to do so, we need to take into account the frequency-dependent absorption factor.  We first consider the WB functional. We assume for concreteness that the noise contribution is larger than the single scattering one, but the same  analysis holds for the reversed situation. Let us fix the SNR at one. Hence, $\lambda_m$ has to be such that $e^{-\Sigma L/2}=\sigma_n$. Writing $\sigma_n=e^{-\tau_n}$, this yields $\Sigma L=2 \tau_n$. We will see in section \ref{mean} that $\Sigma$ admits the following expression 
$$
\Sigma \sim  \frac{ \sigma_0^2 \eta^{3-\delta}}{\eps (1-\frac{\delta}{2})}=\frac{ \sigma_0^2 \ell_c^{3-\delta} L}{\lambda_m^{4-\delta} (1-\frac{\delta}{2})},
\qquad \textrm{which gives} \qquad 
\lambda_m^{4-\delta} \sim \frac{ \sigma_0^2 \ell_c^{3-\delta} L}{\tau_m (1-\frac{\delta}{2})}.
$$
Hence, as can be expected, $\lambda_m$ decreases as the fluctuations of the random medium and the intensity of the noise decrease (i.e. $\sigma_0 \downarrow$ and $\tau_n \uparrow$.) Remark as well that $\lambda_m$ decreases as $\ell_c$ does, and in a faster way compared to $\sigma_0$ or $\tau_n$. This is also expected since the limit $\ell_c \to 0$ corresponds to the homogenization limit in which the wave propagates  in a homogeneous medium provided $\sigma_0$ is small. The measurements are therefore primarily coherent and both the CB and WB functionals perform well. 

Let us consider now the CB functional and let $\delta=0$. Since for $(L/(N_C \lambda))\wedge \mu \gg 1$, $SNR^C$ is greater than $SNR^W$, we may expect in principle to find a lower minimal wavelength. A way to exploit this fact is to consider a central wavelength of the source $\lambda_m / \alpha$, for the $\lambda_m$ above and some $\alpha>1$, and to compute correlations over a sufficiently small domain in order to gain stability, but a domain not too large so that the resolution is still better than $\lambda_m$. The resolution of the CB functional being $L/N_C$, we choose $N_C= L \beta / \lambda_m$, with $\beta>1$ so that $\lambda_m/ \beta$ is smaller than $\lambda_m$. The prefactor in the SNR is then the square of $(L\alpha/ (N_C \lambda_m))\wedge \mu =(\alpha/\beta) \wedge \mu$ that we suppose is equal to $(\alpha/\beta)^2$. We find
$$
SNR^C(0) \sim \sigma_n^{4\alpha-1} (\alpha/\beta)^2,
$$
and we look for $\alpha>1$ and $\beta>1$ such that $\sigma_n^{4\alpha-1} \alpha/\beta=1$. Since $\sigma_n \leq 1$, this is possible only when $\alpha$ is not too large (so that the central wavelength of the source cannot be too small, otherwise absorption is too strong) and when $\sigma_n$ is not too small either. When $\sigma_n$ is below the threshold, then the minimal wavelength is the same for the CB and WB functionals. Hence, compared to the WB functional, the averaging effects of the CB functional can be exploited in order to improve the optimal resolution only when the noise is significant. Finding the optimal $N_C$ is actually a difficult problem, and  is addressed numerically in \cite{Borcea-adapt}.

\paragraph{Effects of long-range correlations.} The strongest such effect is seen in $\Sigma$ since $\Sigma \to \infty$ as $\delta \to 2$. Hence, as the fluctuations get correlated over a larger and larger spatial range, the mean free path decreases and the amplitude of the signal becomes very small. This means that coherent-based functionals are therefore not efficient in random media with long-range correlations, and inversion methodologies based on transport equations \cite{BR-AM, BP-M3AS} that rely on uncoherent information should be used preferrably.

There is also a loss of stability in the variance of the single scattering term for the CB functional, which for sufficiently long correlations (i.e. $\delta>1$) becomes larger than in the short-range case. In such a situation, for the SNR of the CB functional to be larger than that of the WB functional, one needs to decrease the resolution by a factor greater than in the short-range case. Notice moreover that there is an effect of the long-range dependence on the term $(\ell_c/\lambda)^{3-\delta}$ that measures the distance to the homogenization regime. As the medium becomes correlated over larger distances, the variance increases, which suggests that the homogenization regime becomes less accurate.

\section{Imaging functionals and models} \label{meas_imag}
We introduce in sections \ref{WBfunc} and \ref{CBfunc} the WB and CB functionals, and in section \ref{models} our models for the measurements.

\subsection{Expression of the WB functional} \label{WBfunc}
We give here the expression of the WB functional for generic solutions of the wave equation, which we denote by $p$ and its time derivative $\partial_t p$. The solution to the homogeneous wave equation with regular initial conditions $(p(t=0)=q_0, \partial_ t p(t=0)= q_1)$, for $(q_0,q_1)$ given,  reads formally in three dimensions
$$
p(t)=\partial_t G(t,\cdot) * q_0 + G(t,\cdot) * q_1, \qquad G(t,x)=\frac{1}{4 \pi |x|} \delta_0(t-|x|),
$$
where  $*$ denotes convolution in the spatial variables and $\delta_0$ the Dirac measure at zero. In the Fourier space, this reduces to
$$
\calF p(t,k)= \cos |k| t \; \calF q_0(k)+\frac{\sin |k| t}{|k|}\calF q_1(k),
$$
where we define the Fourier transform as $\calF p(k)=\int_{\Rm^3} e^{-i k \cdot x} p(x)dx$.

From the data $(p(T,x), \partial_t p(T,x))$ at a time $t=T$ for $x \in D$, the natural expression of the WB functional in our setting is obtained by backpropagating the measurements, similarly to the time reversal procedure \cite{Fink-Prada-01}. This leads to the functional
$$
I^W_0(x)=[\partial_t G(T,\cdot) * \un_D p(T,\cdot)- G(T,\cdot) * \un_D \partial_t p(T,\cdot)](x),
$$
where $\un_D$ denotes the characteristic function of the detector. Using Fourier transforms, this can be recast as
$$
I^W_0(x)=\calF^{-1}_{k \to x}\left( \cos |k| T \; \calF (\un_D p(T,\cdot))(k)-\frac{\sin |k| T}{|k|}  \calF ( \un_D  \partial_t p(T,\cdot))(k)\right).
$$
This expression is slightly different from the classical Kirchhoff migration functional because of our different measurement setting. It nevertheless performs the same operation of backpropagation. Assuming that all the wavefield is measured, i.e. $D=\Rm^3$, and that $q_1=0$ so $\calF p(t,k)= \cos |k| t \; \calF q_0(k)$, one recovers the initial condition perfectly, i.e. $I^W_0(x)=q_0(x)$. In practice, the entire wavefield is generally not available and diffraction effects limit the resolution. In order to simplify the calculations, we assume without loss of generality that only the pressure $p(T,x)$ is used for this functional. This modifies $I^W_0$ as
\be \label{expK}
I^W(x)=\calF^{-1}_{k \to x}\left( \cos |k| T \; \calF (\un_D p(T,\cdot))(k)\right).
\ee
This significantly reduces the technicalities of our derivations while very little affecting  the reconstructions. Indeed, for localized initial conditions, the value of the maximal peak of the functional $I^W_0$ is divided by a factor two  compared to the full $I^W_0$: if $D=\Rm^3$, and $\calF p(t,k)= \cos |k| t \; \calF q_0(k)$, then
$$
I^W(x)=\calF^{-1}_{k \to x}\left( (\cos |k| T)^2 \calF q_0(k) \right)=\frac{1}{2}\calF^{-1}_{k \to x}\left( \calF q_0(k) \right)+\frac{1}{2}\calF^{-1}_{k \to x}\left( \cos 2 |k| T\calF q_0(k) \right).
$$
The first term above yields $\frac{1}{2}q_0(x)$ while the second one is essentially supported on a sphere of radius $2T$ far away from the source. We will analyse in section \ref{secWB} the expected value and the variance of the functional $I^W$ for random measured wavefields. 

\subsection{Expression of the CB functional} \label{CBfunc}

We define here the CB functional for a measured wavefield $\bu$. This requires us to introduce the Wigner transform of   $\bu$ \cite{LP,GMMP}  and to decompose it into propagating and vortical modes. Since the initial velocity is identically zero, the amplitude of the latter modes remains zero at all times, see \cite{RPK-WM}. The projection onto the propagating mode is done with the eigenvectors of the dispersion matrix $L(k)=k_i D^i$ where $D^i$ was defined in \fref{eq:hypsyst}. These vectors are given by $b_\pm(k)=(\hat k,\pm 1)/\sqrt{2} $, with $\hat k=k/|k|$ and we define in addition the matrices $B_{\pm}=b_\pm \otimes b_\pm$, where $\otimes$ denotes tensor product of vectors. The full matrix-valued Wigner transform of $\bu$ is defined by
$$ 
W^\eps(t,x,k)=\frac{1}{(2\pi)^3} \int_{\Rm^3} e^{i \, k \cdot y}\,\bu(t,x-\frac{\eps y}{2}) \otimes  \bu(t,x+\frac{\eps y}{2}) \, dy.
$$
 The quantity $W^\eps$ is a real-valued matrix.  In our setting, the field $\bu$ is measured at a time $t=T=1$. We compute the correlations over a domain $\eC \subset D$, so that we do not form the full Wigner transform but only a smoothed version of it given by
 $$ 
W^\eps_S(T,x,k)=\frac{1}{(2\pi)^3} \int_{x\pm \frac{\eps y}{2} \in \eC} e^{i \, k \cdot y}\,\bu(T,x-\frac{\eps y}{2}) \otimes  \bu(T,x+\frac{\eps y}{2}) \, dy.
$$
We only consider the propagating mode associated with the positive speed of propagation $c_0$, the other one can be recovered by symmetry. This mode corresponds to the vector $b_+$ and the associated amplitude  is given by
\be \label{aD}
a_S(T,x,k):=\Tr (W^{\eps}_S (T,x,k))^T B_{+}(k),
\ee
where $(W^{\eps}_D)^T$ denotes the matrix transpose of $W^\eps_S$ and $\Tr$ the matrix trace. The CB functional is then defined in our setting by
$$I^C(x)= \int_{D_x} dk \, a_S(T,x+c_0 T\hat{k},k), \qquad D_x=\{ k \in \Rm^3, (x+c_0 T\hat{k}) \in D\}.
$$
Above, we suppose implicitely that if $(x+c_0 T\hat{k}) \in D$, then $\eC$ is such that $(x+c_0 T\hat{k}) \pm \frac{\eps y}{2} \in D$. The functional is slightly different from the classical coherent interferometry functional of \cite{CINT} because of our measurement setting. However, the two functionals qualitatively perform the same operation, that is the backpropagation of field-correlations along rays. The difference resides in how these correlations are calculated: in our configuration, we have access to 3D volumic measurements at a fixed time so that the 3D spatial Wigner transform is available and is the retropropagated quantity; in \cite{CINT}, measurements are performed on a 2D surface and recorded in time so that the 3D spatial Wigner transform is replaced by a 2D spatio-1D temporal Wigner transform. 

For the stability and resolution analysis, it is convenient to recast $I^C$ in terms of the amplitude $a$ associated with the full Wigner transform (i.e. computed for $\eC=\Rm^3$). We then obtain the expression, using the rescaled variables, $c_0=T=1$:
\be \label{filter}
I^C(x)= \int_{D_x} \int_{\Rm^3} dk dq F^x_\eps(k-q) a(1,x+\hat{k},q), \qquad F^x_\eps(k)= \int_{x\pm \eps y/2 \in \eC} e^{i \, k \cdot y} dy.
\ee
Note that $I^C$ is real-valued.

\subsection{Models for the measurements} \label{models}
We introduce in this section our different models for the measurements. We will define a model for the mean of the measurements, and a model for the statistical instabilities. The latter is obtained by using the single scattering approximation (Born approximation). We start with the mean.

\subsubsection{Model for the mean} \label{mean}
\paragraph{The wavefield.} Denoting by $\Box=\partial^2_{t^2}-\Delta$ the d'Alembert operator and by 
\be \label{balK}
p_B(t)=\partial_t G(t,\cdot) * q_0
\ee
the ballistic part associated to an initial condition $q_0$ (with vanishing initial $\partial_t p$), the solution to \fref{waveq} reads:
$$
p(t,x)=p_B(t,x)- \sigma_0 \Box^{-1} \left[V \left(\frac{\cdot}{\eta \eps} \right)\frac{\partial^2 p}{\partial t^2}
 \right](t,x).
 $$
Obtaining an expression for the expectation of $p$ requires the analysis of the term $\E V ((\cdot/(\eta \eps)) \partial^2_{t^2} p$ which is not straighforward and requires a diagrammatic expansion. Rather, we follow the simpler, heuristic, method of \cite{keller64}, which amounts to iterating the above inversion procedure one more time and getting
\begin{align} \label{p2}
&p(t,x)= \\ \nonumber
&p_B(t,x)- \sigma_0 \Box^{-1} \left[V \left(\frac{\cdot}{\eta \eps} \right)\frac{\partial^2 p_B}{\partial t^2} \right](t,x)+
\sigma_0^2 \Box^{-1} \left[V \left(\frac{\cdot}{\eta \eps} \right)\frac{\partial^2}{\partial t^2} \Box^{-1} V \left(\frac{\cdot}{\eta \eps} \right)\frac{\partial^2 p}{\partial t^2}\right](t,x),
\end{align}
and to replace $p$ in the last term by $p_B$. Since $\E\{ p\} \simeq p_B$ at first order in $\sigma_0$, $p_B$ may be  replaced  by $\E\{ p\}$ to obtain a homogenized equation for $\E\{ p\}$. The result is \cite{keller64} that a harmonic wave $\hat p(\omega,x)$ (Fourier transform in time of $p$) is damped exponentially, i.e. for $|x| \neq 0$,
$$
|\E\{\hat p(\omega,x) \}|\leq C e^{-\gamma(\omega) |x|} \qquad \textrm{with} \qquad \gamma(\omega)=\sigma_0^2 \omega^2 \int_{\Rm^3} \frac{1-\cos (2 \omega |x|)}{ 16 \pi |x|^2} R\left(\frac{|x|}{\eps \eta}\right) dx.
$$
The absorption coefficient $\gamma$ is obtained by adapting the setting of \cite{keller64} to ours. Since our initial condition \fref{CI1} is mostly localized around wavenumbers with frequency $|k_0|/\eps$, we find that the waves are absorbed by a factor 
$$
\gamma\left(\frac{|k_0|}{\eps}\right)=\frac{\sigma_0^2 |k_0|^2 \eta}{\eps} \int_{\Rm^3} \frac{1-\cos (2 \eta |k_0| |x|)}{ 16 \pi |x|^2} R(|x|) dx :=\gamma_\eps.
$$
A Taylor expansion then leads to the classical $|k|^4$ dependency of the absorption associated to the Rayleigh scattering:
$$
\gamma_\eps \simeq \frac{\sigma_0^2 |k_0|^4 \eta^3 }{8 \pi \eps} \int_{\Rm^3}R(|x|) dx =  \frac{\sigma_0^2 |k_0|^4 \eta^3 }{8 \pi \eps} \hat{R}(0).
$$
 We implicitely assumed above that $\hat{R}(0)$ was defined, which holds in random media with short-range correlations but not in media with long-range correlations. The latter case is addressed below. We can relate the latter expression for $\gamma_\eps$ to the mean free path $\Sigma^{-1}:=\Sigma^{-1}(|k_0|)$ defined in the next paragraph in \fref{RTE2} by
$$
\Sigma= \frac{\eta^3 \sigma_0^2}{\eps}\frac{\pi |k_0|^{4}}{2 (2\pi)^3} \int_{S^{2}}\hat{R}(\eta(k_0-|k|\hat p)) d \hat p \simeq \frac{\eta^3 \sigma_0^2}{\eps}\frac{\pi |k_0|^{4}}{(2\pi)^2}  \hat{R}(0).
$$
Hence $\Sigma\simeq 2\gamma(|k_0|/\eps)$.

Obtaining such a relation in the case of long-range correlations is slightly more involved.  To do so, we first notice that
$$R\left(\frac{|x|}{\eps \eta}\right)\simeq \frac{(\eps \eta)^{3-\delta} S(0)}{(2 \pi)^3} \int_{\Rm^3} \frac{e^{i k \cdot x}}{|k|^\delta} dk=c_\delta S(0) (\eps \eta)^{3-\delta}  |x|^{\delta-3},$$
where the constant $c_\delta$ is given by $c_\delta=2^{3-\delta} \pi^{3/2} \Gamma(\frac{3-\delta}{2}) (\Gamma(\frac{\delta}{2}))^{-1}(2 \pi)^{-3}$ \cite{gelfand}, 
$\Gamma$ being the Gamma function. The expression of $\gamma_\eps$ is then
$$
\gamma_\eps \simeq c_\delta S(0) \eta^{3-\delta} \sigma_0^2 \eps^{1-\delta} \int_{\Rm^3} \frac{1-\cos\left(\frac{2 |k_0| |x|}{\eps}\right)}{ 16 \pi |x|^{5-\delta}}  dx=\frac{\sigma_0^2 |k_0|^{4-\delta} \eta^{3-\delta}}{2^\delta \pi^{3/2}\eps}  \frac{\Gamma(\frac{3-\delta}{2})}{2\Gamma(\frac{\delta}{2})}\int_{0}^\infty \frac{1-\cos (r)}{r^{3-\delta}}  dr.$$
Recall that $\delta<2$ so that the above integral is well-defined. The inverse of the mean free path is now in the long-range case
\bee
\Sigma&\simeq&  \frac{\eta^{3-\delta} \sigma_0^2}{\eps}\frac{\pi |k_0|^{4-\delta} S(0)}{2 (2\pi)^3 2^{\delta/2}} \int_{S^{2}} \frac{1}{(1-\hat k_0 \cdot \hat p)^{\delta/2}} d \hat p=\frac{\eta^{3-\delta} \sigma_0^2}{\eps}\frac{\pi^2 |k_0|^{4-\delta} S(0)}{ (2\pi)^3 2^{\delta/2}} \int_{-1}^{1} \frac{1}{(1-x)^{\delta/2}} d x\\
&=&\frac{ \sigma_0^2 \eta^{3-\delta} |k_0|^{4-\delta}}{\eps}\frac{S(0)}{  2^{\delta} 4 \pi (1-\frac{\delta}{2})}.
\eee
We used again the fact that $\delta<2$ to make sense of the integral. Since the following relation holds 
$$
\frac{\Gamma(\frac{3-\delta}{2})}{\sqrt{\pi}\Gamma(\frac{\delta}{2})}\int_{0}^\infty \frac{1-\cos (r)}{r^{3-\delta}}  dr=\frac{1}{ 4 (1-\frac{\delta}{2})},
$$
we recover that $\Sigma\simeq 2\gamma(|k_0|/\eps)$ in the long-range case. We will therefore systematically replace $\gamma(|k_0|/\eps)$ by $\Sigma/2$ in the sequel. The expression of $\E \{p\}(t,x)$ is obtained by Fourier-transforming back $\E \{\hat p\}(\omega,x)$ and using the fact that $|x|=c_0 t$ since the ballistic part is supported on the sphere of radius $|x|=c_0 t$. This yields
\bea
\label{meanH}
\E \{p\}(t,x) &\simeq& e^{- c_0 \Sigma  t /2}p_B(t,x)\label{pS}.
\eea
This is our model for the average pressure.
\paragraph{The Wigner transform.} It is shown in \cite{RPK-WM} that the Wigner transform $W^\eps$ of a wavefield $\bu$ satisfies the following system 
\be \label{eq:Wigner}
\frac{\partial W^\eps}{ \partial t}+ \left(\calQ_1+\frac{\calQ_2}{\eps} \right) W^\eps=\left(\calP_1+\frac{\calP_2}{\eps} \right) W^\eps:= \calS^\eps W^\eps,
\ee
with
\bee
\calQ_1 W&=&\frac{1}{2} \left(D^j \pdr{W}{x_j}+\pdr{W}{x_j} D^j\right), \qquad \calQ_2 W=i k_j D^j W-i W k_j D^j\\
\calP_1 W&=&\frac{\sigma_0}{2} \int \int \frac{dy dp e^{i p \cdot y}}{(2 \pi)^d} \left[ \calV\left( \frac{x+\eps y}{\eps \eta}\right)D^j \pdr{W(k+\frac{1}{2}p)}{x_j}\right.\\
&&\left. \hspace{5cm}+\pdr{W(k-\frac{1}{2}p)}{x_j} D^j\calV\left( \frac{x+\eps y}{\eps \eta}\right)\right]\\
\calP_2 W&=&i\sigma_0  \int \int \frac{dy dp e^{i p \cdot y}}{(2 \pi)^d} \left[ \calV\left( \frac{x+\eps y}{\eps \eta}\right)[k+p/2]_jD^j W(k+p/2)\right.\\
&&\left. \hspace{5cm}-W(k-p/2) [k-p/2]_j D^j\calV\left( \frac{x+\eps y}{\eps \eta}\right)\right].
\eee
Noticing that 
$$
\int dy e^{i p \cdot y} \calV\left( \frac{x+\eps y}{\eps \eta}\right)=\eta^3 e^{-i \frac{p \cdot x}{\eps}} \hat{\calV}(-\eta p),
$$
we can recast the r.h.s of (\ref{eq:Wigner}) as 
\bee
\calS^\eps W^\eps &=& F_\eps *_p \left(D^j \left[\frac{\eps}{2}\pdr{W^\eps}{x_j} +i p_j W^\eps\right] \right)+\left(\left[\frac{\eps}{2}\pdr{W^\eps}{x_j} -i p_j W^\eps\right] D^j \right) *_p \overline{F_\eps}\\
F_\eps(x,p)&=&\frac{\sigma_0} {\eps} \left(\frac{\eta} {\pi}\right)^3  e^{\frac{2 i p \cdot x}{\eps}} \hat{\calV}(2\eta p).
\eee
It is then well established \cite{RPK-WM} that a good approximation of $\E \{  W^\eps\}$ is
$$
\E \{  W^\eps\} \simeq \overline{W}_0^\eps:=\sum_{\pm } \bar a_\pm^\eps B_\pm,
$$
where the matrices $B_\pm$ were introduced in section \ref{CBfunc} and the amplitudes $\bar a_\pm^\eps(t,x,k)$ satisfy the following radiative transfer equation
\begin{align} \label{RTE}
&\frac{\partial \bar a_\pm^\eps}{ \partial t} \pm  c_0 \hat k \cdot \nabla_x \bar a_\pm^\eps =Q( \bar a_\pm^\eps),\\ \nonumber
&Q(\bar a_\pm^\eps)(k)=\int_{\Rm^3} \sigma(k,p) (\bar a_\pm^\eps(p)- \bar a_\pm^\eps(k)) \delta_0(c_0|k|-c_0|p|) d p,\\\nonumber
&\sigma(k,p)= \frac{\eta^3 c_0^2 \sigma_0^2}{\eps}\frac{\pi  |k|^2}{2 (2\pi)^3}\hat{R}(\eta(k-p)).
\end{align}
Above, $\delta_0$ is the Dirac measure at zero. The equation (\ref{RTE}) is supplemented with the initial condition
$$
\bar a_\pm(t=0)=a_{0,\pm}^\eps= \Tr (W^{\eps,0})^T B_{\pm},
$$
where $W^{\eps,0}$ is the Wigner transform of the initial wavefield $\bu(t=0)$. Since $\bar a_-^\eps(k)=\bar a_+^\eps(-k)$, we focus only on the mode $\bar{a}_+^\eps$ and drop both the $+$ lower script and the $\eps$ upper script for notational simplicity. We now need to compute the Wigner transform of the initial condition. We consider an approximate expression that simplifies the analysis of the imaging functionals. The scalar Wigner transform of $p_0^\eps$, denoted by $w_0^\eps$, reads:
\bee
w_0^\eps(x,k)&\sim&\frac{|k_0|^4}{(2\pi)^3}\int_{\Rm^3} \int_{S^2 }\int_{S^2 } dy d\hat{k}_1 d\hat{k}_2e^{i (k-|k_0| (\hat{k}_1+\hat{k}_2)/2) \cdot y} e^{i |k_0|x \cdot(\hat k_1 -\hat k_2)/\eps}  \\
&&\hspace{3cm} \times \hat{g}\left( \frac{\hat k_1 \cdot (x-\eps y/2)}{\eps \mu}\right)  \overline{\hat{g}}\left( \frac{\hat k_2 \cdot (x+\eps y/2)}{\eps \mu}\right).
\eee
Since the bandwidth parameter $\mu$ is such that $\mu\eps \gg \eps$, we can separate the scales of the $\hat{g}$ terms and the oscillating exponentials. We can then state that, as is classical for Wigner transforms, different wavevectors lead in a weak sense to negligible contributions because of the highly oscillating term $e^{i |k_0| x \cdot(\hat k_1 -\hat k_2)/\eps}$. The leading term in the initial condition is therefore
\bee
w_0^\eps(x,k)&\simeq& \int_{\Rm^3} \int_{S^2 } dy d\hat k_1  e^{i (k-|k_0| \hat k_1) \cdot y} \hat{g}\left( \frac{\hat k_1 \cdot (x-\eps y/2)}{\eps \mu}\right)  \overline{\hat{g}}\left( \frac{\hat k_1 \cdot (x+\eps y/2)}{\eps \mu}\right)  \\
&:=& \mu^3 \int_{S^2} d \hat k_1 W_{\hat k_1} \left(\frac{x}{\eps \mu},\frac{k-|k_0| \hat k_1}{\mu^{-1}}\right),
\eee
where 
$$
W_{\hat k_1}(x,k) = \int_{\Rm^3} dy e^{i k \cdot y} \hat{g}\left(\hat k_1 \cdot (x- y/2)\right) \overline{\hat{g}}\left(\hat k_1 \cdot (x+y/2)\right).
$$
The effect of the function $w_0^\eps$ in the phase space is essentially to localize the variable $x$ around $\eps \mu$ and the variable $|k|$ on a shell of radius $|k_0|$ with width $\mu^{-1}$. This is what is expected since the initial condition for the wave equation is isotropic, oscillates at a frequency $\eps^{-1} |k_0|$ and has a bandwidth of order $(\eps \mu)^{-1}$. The fact that $\mu \ll 1/\sqrt{\eps}$ implies that the localization is stronger in the spatial variables than in the momentum variables, which can be seen as a broadband property. For the sake of simplicity, we replace the initial condition by a function that shares the same properties but has an easier expression to handle. We write:
$$
a_{0,\pm}^\eps=\mu w_0\left( \frac{x}{\eps \mu},\mu(|k|-|k_0|)\right),
$$
where the function $w_0$ is smooth. With $c_0=1$, we then recast \fref{RTE} as
\begin{align} \label{RTE2}
&\frac{\partial \bar a}{ \partial t} +\hat k \cdot \nabla_x \bar a + \Sigma(k) \bar a=Q_+( \bar a), \qquad \bar{a}_\pm(t=0,x,k)=\mu w_0\left( \frac{x}{\eps \mu},\mu(|k|-|k_0|)\right)\\\nonumber
&Q_+(\bar a)(k)=\int_{\Rm^3} \sigma(k,p) \bar a(p) \delta_0(|k|-|p|) d p, \qquad \Sigma(k)= \frac{\eta^3 \sigma_0^2}{\eps}\frac{\pi |k|^{4}}{2 (2\pi)^3} \int_{S^{2}}\hat{R}(\eta(k-|k|\hat p)) d \hat p.
\end{align}
When $ \sigma_0^2= \eps$, $\eta=1$, we recover the usual equation for the weak coupling regime. Note that since $R(x)=R(|x|)$, we have $\hat{R}(k)=\hat{R}(|k|)$ and $\Sigma(k)=\Sigma(|k|)$. Going back to the measured amplitude $a_D$ defined in (\ref{aD}), an accurate approximation of its mean is therefore
\be
\E\{a_ S\} =\Tr (\E \{W^{\eps}_S\})^T B_{+} =\Tr (F^x_\eps *_k \E \{W^\eps\})^T B_{+}\simeq  F^x_\eps *_k \bar{a}, \label{meana}
\ee
where the filter $F^x_\eps$ was defined in \fref{filter}. The integral solution to  (\ref{RTE2}) reads, denoting by $T_t \bar a(x,k):=\bar a(x-c_0 t \hat k,k)$ the free transport semigroup:
\be \label{integ}
\bar a(t,x,k) =T_t a_0^\eps(x,k) e^{-\Sigma(|k|) t}+  \int_0^t ds e^{-\Sigma(|k|) (t-s)}T_{t-s} Q_+( \bar a)(x,k).
\ee
The average $\bar a(t,x,k)$ is the sum of a ballistic term and the multiple scattering contribution.

\subsubsection{Model for the random fluctuations in the measurements} \label{flucmod}
\paragraph{The wavefield.} We follow the Born approximation, which consists in retaining in \fref{p2} the terms at most linear in $V$. In order to take into account the absorption as explained in section \ref{results}, we replace the ballistic term in \fref{p2} by $\E\{p\}$. In doing so, both the average $\E\{p\}$ and the fluctuations are exponentially decreasing. The random instabilities on the pressure are then generated by the term
\be \label{randH}
\delta p^\eps (t,x)=- \sigma_0 e^{-c_0 \Sigma t/2 }\Box^{-1} \left[V \left(\frac{\cdot}{\eta \eps} \right)\frac{\partial^2 p_B}{\partial t^2} \right](t,x).
\ee
We suppose that the additive noise on the wavefield $\bu$ has the form $\sigma_n \bn^\eps(x)=\sigma_n (\bn_v^\eps,n_p^\eps)(x)=\sigma_n (\bn_v(\frac{x}{\eps}),n_p(\frac{x}{\eps}))$ and is independent of the random medium. For simplicity, we assume that the noise entries are real and have the same correlation structure, $\E\{ \bn^\eps_i(x) \bn^\eps_j(y)\}=\Phi(\frac{x-y}{\eps})$, $i,j=1,\cdots,4$, where $\Phi(x) \equiv \Phi(|x|)$ and is smooth. Using \fref{meanH} and \fref{randH}, our model for the measurements is therefore
\be \label{modH}
p(t,x)=\E \{p\}(t,x) +\delta p^\eps (t,x)+\sigma_n n_p^\eps(x).
\ee

\paragraph{The Wigner transform.} 
Let $a^\eps$ be the projection of the full Wigner transform $W^\eps$ onto the $+$ mode, i.e. $a^\eps=\Tr (W^{\eps})^T B_{+}$, so that, for the $a_S$ defined in \fref{aD}, we have $a_S=F^x_\eps *_k a^\eps$. We already know from \fref{meana} that $\E\{a_S\} \simeq F^x_\eps *_k \bar{a}$, with $\bar a$  the solution to the radiative transfer equation (\ref{RTE2}). We subsequently write $a^\eps=\bar a + \delta a^\eps$, where $\delta a^\eps$ accounts for the random fluctuations. The simplest model for $\delta a^\eps$ is obtained for the single scattering approximation which consists in retaining in $W^\eps$ only terms at most linear in $V$, so that the related variance will be at most linear in $\hat{R}$. 
This leads to defining the random term $W^\eps_1$
\be \label{correc}
\frac{\partial W^\eps_1}{ \partial t}+ \left(\calQ_1+\frac{\calQ_2}{\eps} \right) W^\eps_1= \calS^\eps W^\eps_0,
\ee
with vanishing initial conditions and where
$$
\frac{\partial W^\eps_0}{ \partial t}+ \left(\calQ_1+\frac{\calQ_2}{\eps} \right) W^\eps_0=0, \qquad W^{\eps}_0(t=0)=W^{\eps,0}.
$$
As was the case for the pressure, such a definition of $W_1^\eps$ does not take into account the absorption factor $e^{-c_0 \Sigma t}$. We thus formally correct this and set
\be \label{perturb}
\delta a^\eps= e^{-c_0 \Sigma t}\Tr (W^{\eps}_1)^T B_{+}.
\ee
The fluctuations of the amplitude $\delta a^\eps$ and the wavefield can be related by writing the field $\bu$ as $\E\{\bu\}+ \delta \bu$, where $\delta \bu$ is obtained in the single scattering approximation. The perturbation $\delta a^\eps$ is then the projection on the $+$ mode of the sum of  Wigner transforms $W[\E\{\bu\},\delta \bu]+W[\delta \bu, \E\{\bu\}]$. Taking into account the external noise $\sigma_n \bn^\eps(x)$ and denoting by $a_n^\eps$ the projection of $W[\E\{\bu\}, \bn^\eps ]+W[\bn^\eps, \E\{\bu\}]$ on the $+$ mode as in \fref{perturb}, our complete model for the measurements in the single scattering approximation is therefore
\bea \label{meas}
a^\eps(t,x,k)&=&\bar a(t,x,k)+\delta a^\eps(t,x,k)+ \sigma_n a^\eps_n(t,x,k).
\eea

\section{Analysis of the  CB functional} \label{secCB}
We compute in this section the mean and the variance of the CB functional. 

\subsection{Average of the functional} \label{meancint} Assume for the moment that our measurements have the form, for $g$ a smooth function:
$$
a_S=F^x_\eps *_k g, \qquad g(x,k)=g_0(x-c_0 T \hat k,k):=T_T g_0(x,k).
$$
Above, $F^x_\eps(k)$ is the filter defined in \eqref{filter} and $T_t$ is the free transport semigroup introduced in \fref{integ}. Plugging this expression into the CB functional yields:
\bee
I^C[g](x)&=& \int_{D_x} \int_{\Rm^3} dk dq F^x_\eps(k-q) g_0(x+c_0 T (\hat{k}-\hat{q}),q).
\eee
Let us verify first that when $D=\Rm^3$, we recover that $
I^C[g](x) =(2 \pi)^3 \int_{\Rm^3} dq g_0(x,q),
$
so that if $g_0$ is the scalar Wigner transform of a function $\psi$, we find $I^C[g](x)=(2 \pi)^3 |\psi(x)|^2$ and the reconstruction is then perfect. This is immediate since from the definition of $F_\eps^x(k)$ given in (\ref{filter}) we conclude that $F_\eps^x(k)=(2 \pi)^3 \delta_0(k)$ when $\eC=\Rm^3$. When $\eC$ is finite as is the case in practical situations, one does not recover  $\int_{\Rm^3} dq g(x,q)$ but an approximate version of it limited by the resolution of the functional. This point is addressed below.

\paragraph{Amplitude and resolution.} For simplicity, we suppose that the domain $\eC$ is chosen such that $F_\eps^x$ defined in \eqref{filter} is independent of $x$. We suppose in addition  that $\eC$ is a ball centered at zero of a certain radius. We   parametrize $\eC$ by $\gamma \in [0,1]$ such that if $\pm \frac{\eps y}{2} \in \eC$, then we have $|\eps y| \leq r_0 \eps^{1-\gamma}$ for some $r_0>0$. This means that the diameter of the ball is equal to $$r_0 \eps^{-\gamma}:=N_C$$ where the parameter $N_C$ was introduced in section \ref{results}.  Since the (rescaled) detector has a side $\frac{l}{L} <1$, we have necessarily by construction that $r_0<\frac{l}{L}$. We could  easily accommodate for anisotropic domains $\eC$ with additional technicalities. In order to deal with a regular filter, we  smooth out the characteristic function of the unit disk and replace it by some approximate function $\chi(x)\equiv \chi(|x|)$. Rescaling $y$ as $y \to y r_0 \eps^{-\gamma}$, it comes for the filter $F_\eps$:
$$
F_\eps(k)= \frac{r_0^3}{\eps^{3\gamma}} \int_{\Rm^3} e^{i \, r_0\eps^{-\gamma}k \cdot y} \chi(|y|)dy:=\frac{r_0^3}{\eps^{3\gamma}} F\left( \frac{r_0|k|}{\eps^\gamma}\right).
$$

Recall that our measurements read $a_S=F^x_\eps * a^\eps$, where $a^\eps$ is given by \fref{meas}. The total functional is denoted by $I^C[a^\eps](x)$ and its average is given by $I^C[\bar{a}](x)$. The mean of the measurements $\bar{a}$ is the sum of a ballistic term $a^B(t,x,k) =T_t a_0^\eps(x,k) e^{-\Sigma(|k|) t}$ and a multiple scattering term defined in \fref{integ}. Measurements in a homogeneous medium have the form $T_t a_0^\eps(x,k)$. The CB functional is tailored to backpropagate the data along the characteristics of the transport equation, so as to undo the effects of free transport. Since the multiple scattering term is smoother than the ballistic term, one can expect the inversion operation to produce a signal with a lower amplitude for the multiple scattering part than for the ballistic part. We therefore neglect the multiple scattering in the computation of the average functional. Moreover, the resolution of $I^C$ is mostly limited by the filter $F$ and not by the size of the detector since $\eC$ is included in $D$. As a consequence, the limitations due to $D$ on the resolution are negligible and we replace $D_x$ by $\Rm^3$ in the definition of $I^C$. We then find, 
\bee
 \E\{ I^C[a^\eps]\}(x) &\simeq& I^C[a^B](x)\\
&\simeq & \frac{\mu e^{-\Sigma }r_0^3}{\eps^{3 \gamma}} \int_{D_x} \int_{\Rm^3}  F\left( \frac{r_0|k-q|}{\eps^\gamma}\right)w_0\left(\frac{x+\hat{k}-\hat q}{\eps \mu},\mu(|q|-|k_0|) \right)  dk d|q| d \hat q
\eee
so that the reconstructed image is essentially obtained by a convolution-type relation in the $k$ variable between a localizing kernel and the source. This is clear when $\eps^\gamma \ll \eps \mu$ where 
\bee
 \E\{ I^C[a^\eps]\}(x) &\simeq &  e^{-\Sigma} \int_{\Rm^3} \int_{\Rm^3}  F\left(|k|\right)w_0\left(\frac{x+\frac{\eps^\gamma}{|k_0|r_0} (k-(\hat q \cdot k) \hat q)}{\eps \mu},|q| \right)  dk d|q| d \hat q.
\eee
We can then define the resolution of the functional to be the scale of the localization, which is $\eps^\gamma /r_0=N_C^{-1}$, or $L/N_C$ in dimension variables. Hence, the larger the domain on which the correlations are calculated, the better is the resolution. When the domain is large, more phase information about the wavefield is retrieved, and this is the most useful information to achieve a good resolution. On the contrary, when $\gamma=0$, the phase is lost and imaging is performed mostly using the singularity of the Green's function, which lowers the resolution.

Regarding the amplitude of the signal, we set $x=0$ and compute the value of the peak. If the bandwidth parameter $\mu=\eps^{-\alpha}$ with $\alpha>1-\gamma$, we find 
\bee
 \E\{ I^C[a^\eps]\}(0) &\simeq & e^{-\Sigma}\int_{\Rm^3} \int_{\Rm^3}  F\left( |k|\right)w_0\left(0,|q| \right)  dk dq \simeq  e^{-\Sigma}.
\eee
In this case, the convolution is done at a smaller scale than the spatial support of $w_0$, so that there is no geometrical loss of amplitude with respect to $w_0(x=0)$. When $\alpha<1-\gamma$, convolution is done at a larger scale and we have
\bee
 \E\{ I^C[a^\eps]\}(0) &\simeq & r_0^2 \mu^2 \eps^{2(1-\gamma)}e^{-\Sigma}\int_{\Rm^d} \int_{\Rm^3}  F\left( |k|\right)w_0\left(e_1 \theta_1+e_2 \theta_2,|q| \right)  d|k| d\theta_1 d\theta_2 d|q|\\
&\simeq & r_0^2\mu^2 \eps^{2(1-\gamma)} e^{-\Sigma}.
\eee
Above, $\theta_1$ and $ \theta_2$ are such that $\hat k \simeq \hat q+\eps \mu(e_1 \theta_1+e_2 \theta_2)$ with $e_1 \cdot \hat q=e_2 \cdot \hat q=e_1 \cdot e_2=0$ and $|e_1|=|e_2|=1$. The case $\alpha=1-\gamma$ follows similarly. As a conclusion of this section, we therefore have
\be \label{ampsig}
\E\{ I^C[a^\eps]\}(0) \simeq e^{-\Sigma} \eps^{2(1-\gamma)} \left( (r_0 \mu) \wedge \eps^{\gamma-1}\right)^2.
\ee

\subsection{Variance of the functional}
We compute in this section the variance of the CB functional with our measurements \fref{meas}. For simplicity, we replace $c_0$ and $T$ by their actual values, $c_0=T=1$. The total variance is defined by 
\bee
V^C_{\textrm{tot}}(z)&=& \E \left\{I^C[a^\eps]^2(z) \right\}-(\E \left\{I^C[a^\eps](z) \right\})^2=\E \left\{I^C[\delta a^\eps]^2(z) \right\}+\sigma_n^2 \E \left\{I^C[a^\eps_n]^2(z) \right\}\\
&:=&V^C(z)+V^C_n(z).
\eee
For some function $\varphi$ to be defined later on, let us introduce the following quantity 
$$
w_\eps(t)=\int_{\Rm^{12}} dx_1 dx_2 d q_1 d q_2  J(t,x_1,q_1,x_2,q_2) \varphi(x_1,q_1) \varphi(x_2,q_2),
$$
 where  
$$
J= \E \{\mathcal{W}^\eps\}, \qquad  \mathcal{W}^\eps(t,x_1,q_1,x_2,q_2) =
 \delta a^\eps (t,x_1,q_1)  \delta a ^\eps(t,x_2,q_2).
$$
 The function $J$ is the single scattering approximation of the scintillation function of the Wigner transform and can be seen as a measure of the statistical instabilities. See \cite{B-Ito-04, BLP-JMP, BP-CPDE} for an extensive use of scintillation functions in the analysis of the statistical stability of wave propagation in random media. We have the following technical lemma, proved in section \ref{proofs}, that will be used to obtain an explicit expression for the variance $V^C$:
\begin{lemma} \label{lemscin}$w_\eps$ is given by
$$
w_\eps(t)  = \frac{\sigma_0^2 \eta^3e^{-2 \Sigma}} { 4\pi^3\eps^2}  \int_{\Rm^3} dk  \hat{R}(2\eta k) |H_\eps(t,k)|^2 ,
$$
where the function $H_\eps$ reads
\begin{align*}
&H_\eps(t,k)= \frac{1}{2}\int_0^t \int \int ds dx dq e^{2 i \frac{k\cdot x}{\eps}} (\frakH^\eps_1+\frakH^\eps_2)(t,k,s,x,q)\\
&\frakH^\eps_1=\sum_{\sigma_1,\sigma_2 =\pm 1} \sigma_1 \varphi_{\sigma_2} f_{\sigma_1,\sigma_2}, \qquad 
\frakH^\eps_2=\hat k \cdot \hat{q} \sum_{\sigma_1,\sigma_2 =\pm 1} \varphi_{\sigma_2} f_{\sigma_1,\sigma_2} \\
&f_{\sigma_1,\sigma_2}=\mu \left(\frac{1}{2 \mu} \widehat{q} \cdot \nabla_x +i\sigma_2 |q| \right)w_0\left(\frac{x-\sigma_1 s\hat q}{\mu\eps},\mu(|q|-|k_0|)\right)\\
&\varphi_{\sigma_2} =\varphi(x+(t-s)(\widehat{q+\sigma_2k}),q+\sigma_2 k).
\end{align*}
\end{lemma}
\bigskip

If we set 
\be \label{testf}
\varphi_z(x,q)=\int_{D_z}dk F_\eps(k-q) \delta_0(x-\hat{k}-z)
\ee
in the definition of $w_\eps$, then
\bea \nonumber
V^C(z)&\sim &\int_{\Rm^{12}}dp d q dk dk' F_\eps(k-q) F_\eps(k'-p) \E  \{\delta a ^\eps(1,z+\hat{k},q) \delta a^\eps(1,z+\hat{k}',p)\}\\\label{var1}
&\sim& \int_{\Rm^{12}}  dp d q du dv J(1,u,q,v,p) \varphi_z(u,q) \varphi_z(v,p)=w_\eps(1).
\eea

The analysis of the variance is somewhat technical in that the scintillation is integrated against singular test functions (due to the application of the CB functional) while a classical stability analysis amounts to integrating against regular test functions. As before, we assume that the detector $D$ is large enough so that its effects on the variance can be neglected and we replace it by $\Rm^3$ as a first approximation. We have the following two propositions proved in sections \ref{varproof} and \ref{secvarcintn}:
\begin{proposition} \label{varcint} The variance of the CB functional for the single scattering contribution satisfies
$$
V^C(z)  \sim e^{-2\Sigma}\sigma_0^2 \eta^2 \eps^{4(1-\gamma)} \left(\eta^{1-\delta}\eps^{(1-\delta) \wedge 0} \mu^{-2}((r_0\mu) \wedge \eps^{\gamma-1})^4\right) \vee \left( r_0^4 \eta^{(1-\delta) \wedge 0}\right),\quad |z| \leq \eps \mu 
$$
and $V^C$ is mostly supported on $|z| \leq \mu \eps$.
\end{proposition}
The above result gives the optimal dependency on the parameters $\eta$, $\eps$, $\delta$, $\gamma$ and $\mu$. Regarding the contribution of the noise, we have:
\begin{proposition} \label{varcintn} The variance of the CB functional for the noise contribution satisfies
$$
V^C_n(z) \sim e^{-\Sigma}\sigma_n^2 r_0^4 \eps^{4(1-\gamma)} \quad \textrm{for} \quad |z| \leq \eps \mu,
$$
and $V^C_n$ is mostly supported on $|z| \leq \mu \eps$.
\end{proposition}
As before, the result is optimal and the variance decreases as $\gamma$ goes to zero and resolution is lost. Notice that both variances are essentially supported on the support of the initial source. We now turn to the WB functional.
\section{Analysis of the WB functional} \label{secWB}
In this section we compute  the mean and the variance of the WB functional. 
\subsection{Average of the functional} Using our model \fref{modH} and the definition of the functional \fref{expK}, we have
$$
\E\{I^W [p]\}(z) \simeq e^{-\Sigma /2} \E \{I^W [p_B]\}(z).
$$
The cross-range resolution of $I^W$  is the same as the classical Kirchhoff functional and given by the Rayleigh formula $\lambda L/l$.

\subsection{Variance of the functional} The variance is defined by 
$$
V^W(z)= \E \left\{|I^W[p](z)-\E \left\{I^W[p] \right\}(z)|^2 \right\}.
$$
As before, we distinguish the contribution of the single scattering, denoted by $V^W$, from the one of the noise, denoted by $V^W_n$. We then have the following propositions, proved in sections \ref{proofpropvarK} and \ref{proofpropvarKn}:
\begin{proposition} \label{propvarK} The variance of the WB functional for the single scattering contribution satisfies:
$$
V^W(z) \sim e^{-\Sigma}\sigma_0^2\quad \textrm{for} \quad |z| \leq \eps \mu,
$$
and $V^W$ is mostly supported on $|z| \leq \mu \eps$.

\end{proposition}
As for the CB functional, the results are optimal. Regarding the contribution of the noise, we have:
\begin{proposition} \label{propvarKn} The variance of the WB functional for the noise contribution satisfies for all $z$:
\bee
V^W_n(z) &\sim& \sigma_n^2.
\eee
\end{proposition}
In the latter case, the variance is uniform in $z$ and not just supported on the support of the source. This is due to the fact that, contrary to the CB functional which considers the correlation of the noise with the coherent wavefield (i.e. $V^C_n$ involves $ (p^B)^2 (\sigma_n n_p^\eps)^2$), the variance of the noise for the WB functional does not involve the average wavefield and the information about the source is lost (i.e. $V^W_n$ involves only $ (\sigma_n n_p^\eps)^2$). Remark also that $V_n^W$ and $V_n^C$ are comparable when $\gamma=1$, namely when both functionals have the same resolution. $V_n^C$ decreases when the resolution worsens.

\medskip

The proof of theorem \ref{theo} is then a straightforward consequence of \fref{ampsig}, propositions \ref{varcint}, \ref{varcintn}, \ref{propvarK}, and \ref{propvarKn} after recasting $\eps$, $\eta$ and $r_0 \eps^{-\gamma}$ in terms of $\lambda$, $L$, $\ell_c$ and $N_C$.
\section{Proofs} \label{secproofs}
\subsection{Proof of lemma \ref{lemscin}} \label{proofs}
We start from (\ref{correc}) using the notations of section \ref{models}. We project the initial condition $W_0^\eps$ onto the propagating modes and introduce the related amplitudes $a_\pm$ (recall there are no vortical modes since the initial velocity is zero):
$$
W^\eps_0(t,x,k)=\sum_{\pm} a_{\pm}(t,x,k) B_{\pm}(k), \qquad B_{\pm}=b_\pm \otimes b_\pm \qquad b_\pm=\frac{1}{\sqrt{2}} (\hat k,\pm 1)^T,
$$
with
$$
\frac{\partial a_\pm}{ \partial t} \pm \hat k \cdot \nabla_x a_\pm=0, \qquad a_{\pm}(t=0)= \Tr (W^{\eps,0})^T B_{\pm}=w_0^\mu(x/\eps,|k|-|k_0|),
$$
where $w_0^\mu(x,|k|-|k_0|):=\mu w_0(x/\mu,\mu(|k|-|k_0|))$.  Split then $\calS^\eps$ in \eqref{eq:Wigner} into $\calS^\eps_1+\calS^\eps_2$ with obvious notation. We project \fref{correc} on  the $+$ mode and need to compute $\Tr (\calS^\eps_1 W^{\eps}_0)^T(k) B_{+}(k)$. Direct calculations yield
\begin{align*}
&\Tr (\calS^\eps_1 W^{\eps}_0)^T(k) B_{+}(k)=\frac{1}{4}\sum_{\pm} \int dp f_\eps(x,k-p) (\hat k \cdot \hat p \pm 1)
\left(\frac{\eps}{2} \hat p \cdot \nabla_x a_\pm +i |p| a_\pm \right)\\
&f_\eps(x,p)=\frac{\sigma_0} {\eps} \left(\frac{\eta} {\pi}\right)^3  e^{\frac{2 i p \cdot x}{\eps}} \hat{V}(2\eta p).
\end{align*}
In the same way, we find
$$
\Tr (\calS^\eps_2 W^{\eps}_0)^T(k) B_{+}(k)=\frac{1}{4}\sum_{\pm} \int dp \overline{f_\eps}(x,k-p) (\hat k \cdot \hat p \pm 1) 
\left(\frac{\eps}{2} \hat p \cdot \nabla_x a_\pm -i |p| a_\pm \right)
$$
and finally, recasting $\delta a ^\eps$ by $a_1$ for simplicity of notation:
\begin{align*}
&\frac{\partial a_1}{ \partial t} +\hat k \cdot \nabla_x a_1=A_\eps S_\eps, \qquad S_\eps =\frac{1}{2}\sum_{\pm} (\hat k \cdot \hat p \pm 1) 
\left(\frac{\eps}{2} \hat p \cdot \nabla_x a_\pm +i |p| a_\pm \right),\\
&A_\eps =\Re \int dp f_\eps(x,k-p)S_\eps(x,p).
\end{align*}
Above, $\Re$ stands for real part. The integral equation for $a_1$ reads:
$$
a_1(t,x,k)=\int_{0}^t ds A_\eps S_\eps(t-s,x-s \hat{k},k):=D^{-1}A_\eps S_\eps.
$$
Let us introduce the product function 
$$
  \mathcal{W}^\eps(t,x_1,q_1,x_2,q_2):=a_1(t,x_1,q_1)a_1(t,x_2,q_2).
$$
We then have
\begin{align*}
  \mathcal{W}^\eps(t,x_1,q_1,x_2,q_2) 
  &:= (D^{-1}A_\eps S_\eps)(t,x_1,q_1) (D^{-1}A_\eps S_\eps)(t,x_2,q_2)\\ 
  &= \frac{1}{2}\Re \int_0^t\int_0^t \int \int f_\eps(x_1-s_1\hat{q_1}, q_1-\eta_1)
  f_\eps(x_2-s_2\hat{q_2}, q_2-\eta_2) \\
  &\qquad \times S_\eps(t-s_1, x_1-s_1\hat{q_1}, \eta_1) S_\eps(t-s_2, x_2-s_2\hat{q_2},\eta_2) d\eta_1 d \eta_2 d s_1  d s_2\\
&\quad +\frac{1}{2}\Re \int_0^t\int_0^t \int \int f_\eps(x_1-s_1\hat{q_1}, q_1-\eta_1)
  \overline{f_\eps}(x_2-s_2\hat{q_2}, q_2-\eta_2) \\
  &\qquad\times S_\eps(t-s_1, x_1-s_1\hat{q_1}, \eta_1)   \overline{S_\eps}(t-s_2, x_2-s_2\hat{q_2},\eta_2) d\eta_1 d \eta_2 d s_1  d s_2\\
&:=T_1+T_2.
\end{align*}
Using the fact that $\E \{\hat{V}(\xi)\hat {V}(\nu) \}=  (2 \pi)^{3}\hat{R}(\xi) \delta_0(\xi+\nu),$
we find
\bee
\E \{f_\eps(x,p)  f_\eps(y,q)\}&=&\frac{\sigma_0^2 \eta^3} {\pi^3 \eps^2}   e^{\frac{2 i p \cdot (x-y)}{\eps}} \hat{R}(2\eta p) \delta_0(p+q)\\
\E \{f_\eps(x,p)  \overline{f_\eps}(y,q)\}&=&\frac{\sigma_0 \eta^3} {\pi^3 \eps^2}   e^{\frac{2 i p \cdot (x-y)}{\eps}} \hat{R}(2\eta p) \delta_0(p-q).
\eee
Therefore
\begin{align*}
T_1
    &= \frac{\sigma_0^2 \eta^3} {2 \pi^3 \eps^2} \Re \int_0^t\int_0^t \int  \hat{R}(2\eta (q_1-\eta_1)) e^{\frac{2 i (q_1-\eta_1) \cdot (x_1-s_1 \hat{q}_1-x_2+s_2 \hat{q}_2)}{\eps}}\\
  &\;\times S_\eps(t-s_1, x_1-s_1\hat{q_1}, \eta_1) S_\eps(t-s_2, x_2-s_2\hat{q_2},q_2+q_1-\eta_1) d\eta_1 d s_1  d s_2
\end{align*}
and
\begin{align*}
T_2
    &= \frac{\sigma_0^2 \eta^3} {2 \pi^3\eps^2} \Re \int_0^t\int_0^t \int  \hat{R}(2\eta (q_1-\eta_1)) e^{\frac{2 i (q_1-\eta_1) \cdot (x_1-s_1 \hat{q}_1-x_2+s_2 \hat{q}_2)}{\eps}}\\
  &\;\times S_\eps(t-s_1, x_1-s_1\hat{q}_1, \eta_1)  \overline{S_\eps}(t-s_2, x_2-s_2\hat{q}_2,q_2-q_1+\eta_1) d\eta_1 d s_1  d s_2.
\end{align*}
Hence
\begin{align*}
&  \E\{\mathcal{W}^\eps\}(t,x_1,q_1,x_2,q_2) \\
  & =\frac{\sigma_0^2 \eta^3} {2 \eps^2} \Re \int_0^t\int_0^t \int  \hat{R}(2\eta (q_1-\eta_1)) e^{\frac{2 i (q_1-\eta_1) \cdot (x_1-s_1 \hat{q_1}-x_2+s_2 \hat{q_2})}{\eps}}S_\eps(t-s_1, x_1-s_1\hat{q_1}, \eta_1)\\
&\times \left[S_\eps(t-s_2, x_2-s_2\hat{q_2},q_2+q_1-\eta_1)+\overline{S_\eps}(t-s_2, x_2-s_2\hat{q_2},q_2-q_1+\eta_1) \right]\\
&= \frac{\sigma_0^2 \eta^d} {4 \eps^2}\int_0^t\int_0^t \int  \hat{R}(2\eta k) e^{\frac{2 i k\cdot (x_1-s_1 \hat{q_1}-x_2+s_2 \hat{q_2})}{\eps}}\\
&\times \left[S_\eps(t-s_1, x_1-s_1\hat{q_1}, q_1-k)+\overline{S_\eps}(t-s_1, x_1-s_1\hat{q_1}, q_1+k) \right]\\
&\times \left[S_\eps(t-s_2, x_2-s_2\hat{q_2},q_2+k)+\overline{S_\eps}(t-s_2, x_2-s_2\hat{q_2},q_2-k) \right].
\end{align*}
Accounting finally for the absorption $e^{-\Sigma}$, this allows to obtain the following expression for $w_\eps$: 
$$
w_\eps(t) = \frac{\sigma_0^2 \eta^3 e^{-2\Sigma}} {4 \pi^3\eps^2} \int_{\Rm^3} dk  \hat{R}(2\eta k) |H_\eps(t,k)|^2,
$$
where
\begin{align*}
&H_\eps(t,k)=\int_0^t \int_{\Rm^3}\int_{\Rm^3} ds dx dq e^{\frac{2 i k\cdot x}{\eps}} \left[S_\eps(s, x, q-k)+\overline{S_\eps}(s, x, q+k) \right] \varphi(x+(t-s) \hat{q},q).
\end{align*}
We have finally
\bee
S_\eps (s, x, p)&=&\frac{1}{2}\sum_{\pm} (\hat k \cdot \hat p \pm 1) 
\frac{1}{2} \hat p \cdot (\nabla_x w_0^\mu)\left(\frac{x\mp s\hat p}{\eps},|p|-|k_0|\right)\\
&& +\frac{i}{2}\sum_{\pm} (\hat k \cdot \hat p \pm 1)  |p| w_0^\mu\left(\frac{x\mp s\hat p}{\eps},|p|-|k_0|\right),
\eee
which replaced in $H_\eps$ yields the expression given in the lemma. This is the final result and ends the proof.

\subsection{Proof of proposition \ref{varcint}} \label{varproof} We compute formally the leading term in an asymptotic expansion of the variance, which will give us the optimal dependency in the main parameters. The proof involves the analysis of oscillating integrals coupled to localizing terms. We start from \fref{var1}. We treat only the case $\gamma \in(0,1)$, the cases $\gamma=0$ and $\gamma=1$ follow by using the same approach. Plugging \fref{testf} in \fref{var1}, and assuming the detector is sufficiently large so that we can replace $D_z$ by $\Rm^3$ in a first approximation, we find the following expression for $H_\eps$:
\bee
&&H_\eps(t=1,k)=H_\eps^1(k)+H_\eps^2(k)\\
&&H_\eps^1(k)= \frac{1}{2}\sum_{\sigma_1,\sigma_2=\pm 1} \int_0^{1} \int  ds dq du   F_\eps(u-q) \hat k \cdot (\widehat{ q-\sigma_2k}) e^{2 i k\cdot (\sigma_1 s \widehat{q-\sigma_2k}+\Phi_\eps))/\eps} f_{\sigma_2}\left(\frac{\Phi_\eps}{\eps \mu },q\right) \\
&&H_\eps^2(k)= \frac{1}{2}\sum_{\sigma_1,\sigma_2=\pm 1} \sigma_1 \int_0^{1} \int  ds dq du  F_\eps(u-q) e^{2 i k\cdot (\sigma_1 s \widehat{q-\sigma_2k}+\Phi_\eps))/\eps} f_{\sigma_2}\left(\frac{\Phi_\eps}{\eps \mu},q\right)
\eee
where
\begin{align*}
&\Phi_\eps=z-\sigma_1 s(\widehat{q-\sigma_2k})+s\hat{q}-\hat{q}+\hat{u}\\
&f_{\sigma_2}(x,q) = \left(\frac{1}{2} \widehat{q-\sigma_2k} \cdot \nabla_x +i\mu \sigma_2 |q-\sigma_2k| \right)w_0\big(x,\mu(|q-\sigma_2k|-|k_0|)\big).
\end{align*}
The dependency of the function $f_{\sigma_2}$ on $k$ and $\mu$ is not explicited in order to simplify already heavy notations. Expanding $|H_\eps^1+H_\eps^2|^2$ leads to the following decomposition of $w_\eps$: we can write $w_\eps=w_\eps^1+w_\eps^2+w_\eps^{12}$, where
\bee
&&w^1_\eps(t=1)  = \frac{\sigma_0^2 \eta^3e^{-2\Sigma}} { 4\pi^3\eps^2}  \int_{\Rm^3} dk  \hat{R}(2\eta k) |H^1_\eps(k)|^2
\eee
and $w_\eps^1$ and $w_\eps^{12}$ are obtained in the same fashion. We will focus only on the most technical term $w_\eps^1$ since $w_\eps^2+w_\eps^{12}$ lead either to the same leading order in terms of $\eps$, $\gamma$, $\eta$ and $\sigma_0$ or to a negligible contribution. The term $w_\eps^1$ can itself be decomposed as 
\bee
w_\eps^1&=&\int_0^\eps \int_0^\eps ds_1 d s_2 \;(\cdots)+\int_\eps^1 \int_\eps^1 ds_1 ds_2 \; (\cdots)+\int_\eps^1 \int_0^\eps ds_1 ds_2 \; (\cdots)\\
&:=&v_1^\eps+v_2^\eps+v_3^\eps.
\eee
The term $v_3^\eps$ can be shown to be negligible compared to the first two due to the mismatch of the integration domains, and we thus focus on  $v_1^\eps$ and $v_2^\eps$. We start with the most difficult term $v_2^\eps$. We give a detailed derivation for the case $\mu=\eps^{-\alpha}$ with $\alpha<1-\gamma$ (so that $\eps \mu \ll  \eps^\gamma$) and only state the result for the case $\alpha\geq 1-\gamma$. Let us start by denoting by $(\theta_u^1,\theta_u^2)$ and $(\theta_q^1,\theta_q^2)$ the angles defining $\hat u$ and $\hat q$ with the convention $(\theta_u^1,\theta_u^2) \in (0, 2\pi) \times (0,\pi)$, so that $\hat q=(\cos \theta_q^1 \sin \theta_q^2,\sin \theta_q^1 \sin \theta_q^2,\cos \theta_q^2 )$. After performing the changes of variables $(\theta_u^1,\theta_u^2) \to (\theta_q^1,\theta_q^2)+\eps \mu (\theta_u^1,\theta_u^2)$, we first find that
\be \label{change}
\hat u \simeq\hat q+\eps\mu e_ {1} \theta^1_u + \eps\mu e_2 \theta^2_u \qquad \textrm{with} \qquad e_1 \cdot \hat q=e_2 \cdot \hat q=e_1 \cdot e_2=0, \qquad |e_1|=|e_2|=1.
\ee
Setting in addition $s_1 \to \eps s_1$, $s_2 \to \eps s_2$, $|u| \to |q| + \eps^\gamma |u|/r_0$, $k \to \eps k$, $\xi_1=s_1 \hat q_1$ and $\xi_2=s_2 \hat q_2$, we find for the leading term, using the assumptions $\eps \mu \ll  \eps^\gamma$, $ \eps\mu^2 \ll 1$ and $\mu \gg 1$:
\bea \nonumber
v^\eps_2&\sim& \sigma_0^2 \eta^{3-\delta}\eps^{4(1-\gamma)+1-\delta} \mu^4 r_0^4 e^{-2\Sigma}\\
\nonumber&&\times \int_{\Rm^3} \int_{\eps \leq |\xi_1| \leq 1}\int_{\eps \leq |\xi_2| \leq 1}\frac{dk d \xi_1 d\xi_2}{|k|^\delta |\xi_1|^2 |\xi_2|^2}  (\hat k \cdot \hat{\xi}_2)(\hat k \cdot \hat{\xi}_2) e^{2i k\cdot \{\xi_1-\xi_2\}} \\
&&\hspace{4cm} \times G_{z}\left(\frac{|\xi_1|}{\eps \mu},\hat \xi_1 \right)\overline{G_{z}}\left(\frac{|\xi_2|}{\eps \mu},\hat \xi_2 \right)  \label{v2}
\eea
where
\begin{align*}
&G_{z}(|\xi|,\hat \xi)=S(0)\sum_{\sigma_1,\sigma_2=\pm 1}  \int dY |q|^4 F(|u| \hat \xi ) f^0_{\sigma_2} (z_{\eps\mu}+(1-\sigma_1)\xi+e_ {1} \theta^1_u +  e_2 \theta^2_u,|q|, \hat \xi) \\
&\int dY =\int_{\Rm_+} \int_{\Rm}\int_{\Rm} \int_{\Rm_+} d|q| d\theta^1_u d\theta^2_u d|u|, \qquad e_1 \cdot \hat \xi=e_2 \cdot \hat \xi=e_1 \cdot e_2=0\\
&f^0_{\sigma_2} (x,|q|, \hat \xi)=\left(\frac{1}{2} \widehat{\xi} \cdot \nabla_x +i\mu \sigma_2 |q| \right)w_0\big(x,\mu(|q|-|k_0|)\big),
\end{align*}
and $z_{\eps \mu}=z/\eps\mu$. We will use the following relation:
\be \label{delt}
\int_{\Rm^3} dk\frac{e^{i k\cdot x}}{|k|^\beta}=
\left\{
\begin{array}{l}
(2 \pi)^3 \delta_0(x) \qquad \textrm{if} \quad \beta=0,\\
C_\beta|x|^{\beta-3}  \qquad \textrm{if} \quad 0<\beta \qquad \textrm{and} \qquad \beta \neq 3,5,\cdots,
\end{array}
\right.
\ee
where $\delta_0$ is the Dirac measure at zero and $C_\beta=2^{3-\beta} \pi^{\frac{3}{2}} \Gamma(\frac{3-\beta}{2})(\Gamma(\frac{\beta}{2}))^{-1}$ ($\Gamma$ being the gamma function). Note that \fref{delt} still holds for the case $\beta=3,5,\cdots$ provided the constant $C_\beta$ is adjusted, see \cite{gelfand}. We will not explicit this constant since only the dependency in $x$ matters to us. Hence, we deduce from \fref{delt} that
\be \label{delt2}
\int_{\Rm^3} dk\frac{e^{i k\cdot x}}{|k|^\delta}\hat k_j \hat k_l =
C_\delta' \partial^2_{x_j x_l} |x|^{\delta-1}  \qquad \textrm{for} \quad \delta \in [0,2],
\ee
where $C_\delta'$ is a constant. When $\delta \in (0,1)$, setting $\xi_1 \to \eps \xi_1$ and $\xi_2 \to \eps \xi_2$ in \fref{v2}, it comes using \fref{delt2} and $\mu \gg 1$:
\bee
v^\eps_2&\sim& \sigma_0^2 \eta^{3-\delta}\eps^{4(1-\gamma)} \mu^4 r_0^4 e^{-2\Sigma}\\
&&\times  \int_{1\leq |\xi_1|}\int_{1 \leq |\xi_2|}\frac{d \xi_1 d\xi_2 (\hat \xi_1 \cdot (\xi_1-\xi_2))(\hat \xi_2 \cdot (\xi_1-\xi_2))}{|\xi_1-\xi_2|^{5-\delta} |\xi_1|^2 |\xi_2|^2}  G_{z}(0,\hat \xi_1)\overline{G_{z}}(0,\hat \xi_2)
\eee
plus a term of same order. The integral above is easily shown to be finite using the Hardy-Littlewood-Sobolev inequality \cite{RS-80-2}. A close look at $G_z$ shows that only the term proportional to $\hat \xi \cdot \nabla_x$ is left because of the sum over $\sigma_2$. Setting finally $|q| \to |k_0|+|q|/\mu$ yields $v^\eps_2\sim \sigma_0^2 \eta^{3-\delta}\eps^{4(1-\gamma)} \mu^2 r_0^4 e^{-2\Sigma}$ for  $|z|\leq \mu \eps$. Because of the term $f_{\sigma_2}^0$, it is not difficult to see that $v^\eps_2$ is mostly supported on $|z|\leq \mu \eps$. When $\delta \in [1,2)$ we use \fref{delt2} and directly send $\eps$ to zero in \fref{v2}. The leading term in $G_{z}$ is obtained for $\sigma_1=1$, and we find
$$
v^\eps_2\sim \sigma_0^2 \eta^{3-\delta}\eps^{4(1-\gamma)+1-\delta} \mu^4 r_0^4 e^{-2\Sigma}\int_{|\xi_1| \leq 1}\int_{|\xi_2| \leq 1}\frac{d \xi_1 d\xi_2}{|\xi_1-\xi_2|^{3-\delta} |\xi_1|^2 |\xi_2|^2}  G^0_{z}(\hat \xi_1)\overline{G^0_{z}}(\hat \xi_2),
$$
where
$$
G_{z}^0(\hat \xi)=
S(0)\sum_{\sigma_2=\pm 1} \int dY |q|^4 F(|u| \hat \xi ) f^0_{\sigma_2} (z_{\eps \mu}+e_ {1} \theta^1_u +  e_2 \theta^2_u,|q|, \hat \xi).
$$
Again, summation over $\sigma_2$ in $G^0_{z}$ imply that $v^\eps_2\sim \sigma_0^2 \eta^{3-\delta}\eps^{4(1-\gamma)+1-\delta} \mu^2 r_0^4 e^{-2\Sigma}$ for $|z|\leq \mu \eps$. The case $\delta=0$ gives $v^\eps_2\sim \sigma_0^2 \eta^{3}\eps^{4(1-\gamma)} \mu^2 r_0^4$, and $v_2^\eps(z)$ is mostly supported on $|z| \leq \eps \mu$. Gathering all the previous results, we find for the case $\mu=\eps^{-\alpha}$ with $\alpha<1-\gamma$:
$$
v_2^\eps(z) \sim \sigma_0^2 \eta^{3-\delta}\eps^{4(1-\gamma)+(1-\delta) \wedge 0} \mu^2 r_0^4 e^{-2\Sigma}, \quad |z| \leq \eps \mu,
$$
and $ v_2^\eps(z)$  is essentially supported on $|z| \leq \eps \mu$. When $\alpha\geq 1-\gamma$ (so that $\eps \mu \geq \eps^\gamma$), very similar calculations show that for $|z| \leq \eps \mu$, we have $w_\eps(z) \sim \sigma_0^2\eta^{3-\delta}\eps^{(1-\delta) \wedge 0}\mu^{-2}e^{-2\Sigma} $, and that $v_2^\eps(z)\ll v_2^\eps(0)$ when $|z| \gg \eps \mu$. Combining the last two results, it finally follows that
$$
v^\eps_2(z) \sim \sigma_0^2\eta^{3-\delta} \eps^{4(1-\gamma)+(1-\delta) \wedge 0} \mu^{-2} ((r_0\mu) \wedge \eps^{\gamma-1})^4e^{-2\Sigma},\quad |z| \leq \eps \mu,
$$
and $ v_2^\eps(z)$ is mostly supported on  $ |z| \leq \eps \mu$. 

We consider now the term $v_1^\eps$. After the changes of variables $s_1 \to \eps s_1$ and $s_2 \to \eps s_2$, we find 
\bee
&&v_1^\eps \sim  \sigma_0^2 \eta^3 e^{-2\Sigma}\int_{\Rm^3} dk  \hat{R}(2\eta k) |V^1_\eps(k)|^2
\eee
where
\bee
V_\eps^1(k)= \frac{1}{2}\sum_{\sigma_1,\sigma_2=\pm 1} \int_0^{1} \int  ds dq du   F_\eps(u-q) \hat k \cdot (\widehat{ q-\sigma_2k}) e^{2 i k\cdot (s\hat{q}+(z+\hat{u}-\hat{q})/\eps)} f_{\sigma_2}\left(\frac{z+\hat{u}-\hat{q}}{\eps \mu },q\right) \\
\eee
and $f_{\sigma_2}$ is the same as before. Both cases $\alpha<1-\gamma$ and $\alpha \geq 1-\gamma$ lead to the same order, and we therefore only detail the case $\alpha\geq 1-\gamma$. The change of variables $u \to q+\eps^\gamma u/r_0 $ leads to:
\begin{align*}
&V_\eps^1(k) \sim \sum_{\sigma_1,\sigma_2=\pm 1} \int_0^{1} \int  ds dq \hat F\big(2\eps^{\gamma-1} ((k \cdot \hat q) \hat q-k)/(r_0|q|))\big) \hat k \cdot (\widehat{ q-\sigma_2k}) e^{2 i s k\cdot\hat{q}} f_{\sigma_2}\left(z_{\eps \mu},q\right),
\end{align*}
where $\hat F$ is the Fourier transform of $F$. Denoting then as before by $(\theta_k^1,\theta_k^2)$ and $(\theta_q^1,\theta_q^2)$ the angles defining $\hat k$ and $\hat q$, we perform the change of variables $(\theta_q^1,\theta_q^2) \to (\theta_k^1,\theta_k^2)+r_0 \eps^{1-\gamma} (\theta_q^1,\theta_q^2)$, which yields
$$
\hat q \simeq \hat k+\eps^{1-\gamma}r_0 e_ {1} \theta^1_q + \eps^{1-\gamma}r_0 e_2 \theta^2_q \qquad \textrm{with} \qquad e_1 \cdot \hat k=e_2 \cdot \hat k=e_1 \cdot e_2=0, \qquad |e_1|=|e_2|=1.
$$
 The term $V_\eps^1$ then becomes
\begin{align*}
&V_\eps^1(k) \sim r_0^2 \eps^{2(1-\gamma)} \sum_{\sigma_1,\sigma_2=\pm 1} \int_0^{1} \int  ds |q|^2 d|q|  d\theta_q^1 d\theta_q^2 \;\hat F\big(2 |k|(e_ {1} \theta^1_q +  e_2 \theta^2_q)/|q|\big) \\
& \hspace{6cm}\times \textrm{sign}(|q|-\sigma_2 |k|) e^{2 i s |k|} f_{\sigma_2}\left(z_{\eps \mu},\hat k |q|\right)
\end{align*}
We need now to consider separately $\sigma_2=1$ and $\sigma_2=-1$ in the expression above. Let us start with $\sigma_2=1$, so that  $f_{\sigma_2=1}$ involves 
the term $w_0(z_{\eps \mu},\mu (||q|-|k||-|k_0|))$, which suggests the change of variables $|q| \to |k|+|k_0| + |q|/ \mu$. This brings a factor $\mu^{-1}$ than cancels with the $\mu$ in the definition of $f_{\sigma_2}$. It remains to make sense of the integral in $k$ in the definition of $v_1^\eps$. Finite values pose no problems, and lead, after taking the square, to a term of overall order $\sigma_0^2 r_0^4 \eta^{(3-\delta)} \eps^{4(1-\gamma)}e^{-2\Sigma}$. We focus therefore on large values of $k$. The integral of the term depending on $s$  behaves as $|k|^{-1}$, while the one in $\hat F$ is essentially independent of $|k|$ when $k$. Hence, taking the square, we find a term proportional to
$$
\sigma_0^2 r_0^4 \eta^3 \eps^{4(1-\gamma)}e^{-2\Sigma}\int_{|k| \geq 1} dk  \hat{R}(2\eta k) |k|^{-2} \sim \sigma_0^2 r_0^4 \eta^3 \eps^{4(1-\gamma)}e^{-2\Sigma}\int_{|k| \geq  1} d|k|  S(2\eta k) |\eta k|^{-\delta}.
$$
The integral is finite when $\delta>1$ and leads to a factor $\eta^{3-\delta}$. When $\delta \leq 1$, we set $k \to k /\eta$, which leads to a factor $\eta^2$. The term associated to $\sigma_2=1$ is therefore of order $\sigma_0^2 r_0^4\eta^{(3-\delta) \wedge 2} \eps^{4(1-\gamma)}e^{-2\Sigma}$. Consider now the case $\sigma_2=-1$, so that  $f_{\sigma_2=-1}$ involves now $w_0(z_{\eps \mu},\mu (||q|+|k||-|k_0|))$. When $|k|>|k_0|$, the second argument in $w_0$ never vanishes, which leads to high order terms in $\mu^{-1}$. The leading contribution is thus obtained for $|k| \leq |k_0|$. As mentioned above, bounded values of $k$ leads to an overall order of $\sigma_0^2 r_0^4\eta^{(3-\delta)} \eps^{4(1-\gamma)}e^{-2\Sigma}$, which is negligible compared to the case $\sigma_2=1$ when $\eta \ll 1$, or to the same order when $\eta=1$. It follows therefore that $v_1^\eps$ is  of order $\sigma_0^2 r_0^4\eta^{(3-\delta) \wedge 2} \eps^{4(1-\gamma)}e^{-2\Sigma}$, and with the same arguments as $v_2^\eps$, that it is mostly supported on $|z| \leq \eps \mu$. Owing this, and taking the maximum with $v_2^\eps$ gives the result of the proposition.

\subsection{Proof of proposition \ref{varcintn}} \label{secvarcintn}
We use the same notation as in section \ref{flucmod} and only treat the case $\gamma \in (0,1)$ for brevity. We have  
\bee V^C_n(z)&=&\sigma_n^2 \E \left\{I^C[a^\eps_n]^2(z) \right\}\\
&\sim& \sigma_n^2 \int_{\Rm^{12}}dp d q dk dk' F_\eps(k-q) F_\eps(k'-p) \E  \{a_n^\eps(1,z+\hat{k},q) a_n^\eps(1,z+\hat{k}',p)\}.\eee
Moreover
\bee 
a_n^\eps(1,x,q)&=&\Tr (W[\E\{u\},\bn^\eps]+W[\bn^\eps, \E\{u\}])^T(1,x,q) B^+(q)\\
&=& w[b_+ \cdot \E\{u\}, b_+ \cdot\bn^\eps ]+w[b_+ \cdot\bn^\eps,b_+ \cdot \E\{u\}]\\ 
&:=&\frac{1}{2}(w[\E\{p\}, \bn_p^\eps ]+w[\bn^\eps_p,\E\{p\}])+ w_\bv,
\eee
where $w[\cdot,\cdot]$ denotes the scalar Wigner transform and $w_\bv$ contains the remaining terms not included in the first one. We will not analyse the contributions to the variance of the terms involving $w_\bv$ since they have the same structure as the one involving $w[\E\{p\}, \bn_p^\eps ]+w[\bn^\eps_p,\E\{p\}]$ and yield similar results. Then,
\begin{align}\nonumber
&\E  \{a_n^\eps(1,z,q) a_n^\eps(1,z',p)\}\\\nonumber
&=\frac{1}{4 }\E  \left\{(w[\E\{p\}, \bn_p^\eps ]+w[\bn^\eps_p,\E\{p\}])(z,q)(w[\E\{p\}, \bn_p^\eps ]+w[\bn^\eps_p,\E\{p\}])(z',p) \right\}+R_\eps\\\label{Vnp}
&=\frac{1}{4 }\E  \{w[\E\{p\}, \bn_p^\eps ](z,p)w[\E\{p\}, \bn_p^\eps ](z',q)\} +R'_\eps.
\end{align}
Again, we do not treat the term $R'_\eps$ since it has the same structure as the first term above. We denote by $V^C_{n,p}$ the contribution to the variance of the first term of \fref{Vnp}. It reads, with the notation $u=z+\hat{k}$ and $u'=z+\hat{k}'$:
\begin{align*}
&V^C_{n,p} \sim \sigma_n^2 e^{-\Sigma}\int_{\Rm^{18}}dp d q dk dk' dy dy'F_\eps(k-q) F_\eps(k'-p) e^{i q \cdot y}e^{i p \cdot y'}\\
&\hspace{4cm}p_B(1,u-\frac{\eps}{2} y)p_B(1,u'-\frac{\eps}{2} y')\Phi\left(\frac{u-u'}{\eps}+\frac{1}{2}(y-y')\right).
\end{align*}
In Fourier variables and after setting $\xi \to \eps^{-1}\xi$, $\xi' \to \eps^{-1}\xi'$ and writing the cosines as sum of exponentials, this translates into:
 
\begin{align*}
&V^C_{n,p} \sim \sigma_n^2 e^{-\Sigma}\sum_{\sigma_1, \sigma_2=\pm 1} \int_{\Rm^{15}}d q dk dk' d\xi d\xi' F_\eps(k-q) F_\eps(k'+q-(\xi+\xi')/2) e^{i \frac{u \cdot \xi}{\eps}}e^{i\frac{u' \cdot \xi'}{\eps}}\\
&\hspace{4cm}e^{i \frac{u-u'}{\eps}\cdot (\xi - 2q)} e^{i \frac{\sigma_1|\xi|}{\eps}} e^{i\frac{\sigma_2 |\xi'|}{\eps}}\hat \Phi\left(\xi - 2q\right) g_\mu(|\xi|-|k_0|)g_\mu(|\xi'|-|k_0|).
\end{align*}
Above, we used the fact that $\hat p_B(t,\xi)= \eps^3 g_\mu(\eps(|\xi|-|k_0|/\eps)) \cos t |\xi| $ with $g_\mu(x)=\mu g(\mu x)$. We then decompose $\xi$ and $\xi'$ as $\xi=\xi_\pr \hat k+ \xi_\bot$, $\xi'=\xi_\pr' \hat k+ \xi_\bot'$ with $\xi_\bot \cdot k=\xi_\bot' \cdot k=0$. Denoting by $(\theta_k^1,\theta_k^2)$ and $(\theta_{k'}^1,\theta_{k'}^2)$ the angles defining $\hat k$ and $\hat k'$ and  performing the changes of variables $(\theta_{k'}^1,\theta_{k'}^2) \to (\theta_{k}^1,\theta_{k}^2)+\eps (\theta_{k'}^1,\theta_{k'}^2)$, we find $\hat k'\simeq  \hat k+ \eps e_1 \theta_{k'}^1+\eps e_2 \theta_{k'}^2$ with $e_1 \cdot \hat k=e_2 \cdot \hat k=0$ and $|e_1|=|e_2|=1$. Together with $\xi_\bot \to \sqrt{\eps} \xi_\bot$, $\xi'_\bot\to \sqrt{\eps} \xi_\bot'$, as well as $|\xi_\pr \hat k+ \sqrt{\eps}\xi_\bot| \simeq |\xi_\pr| (1+\eps |\xi_\bot|^2/(2 |\xi_\pr|))$, it comes, since $\mu^2 \ll \eps^{-1} $ and $\eps \ll \eps^\gamma$ as $\gamma \in (0,1)$:
\begin{align*}
&V^C_{n,p} \sim \sigma_n^2 \eps^{4} e^{-\Sigma}\sum_{\sigma_1, \sigma_2=\pm 1} \int dX  F_\eps(k-q) F_\eps(\hat k|k'|+q-\hat k (|\xi_\pr|+|\xi'_\pr|)/2) \\
&\hspace{4cm}\times e^{i(\xi_\pr+\sigma_1 |\xi_\pr|)/\eps}e^{i (\xi_\pr'+\sigma_2 |\xi_\pr'|)/\eps}e^{i \sigma_1 |\xi_\bot|^2/(2|\xi_\pr|)}e^{i \sigma_2 |\xi'_\bot|^2/(2|\xi'_\pr|)}  \\
&\hspace{4cm}\times e^{-i (e_1 \theta^1_{k'}+ e_2 \theta^2_{k'}) \cdot (\xi_\pr \hat k-2q)} e^{i z \cdot (\xi_\eps-\xi'_\eps)/\eps} \\
&\hspace{4cm}\times \hat \Phi(\xi_\pr \hat k- 2q) g_\mu(|\xi_\pr|-|k_0|)g_\mu(|\xi'_\pr|-|k_0|)\\
&\int dX=\int_{\Rm^3}\int_{\Rm^3}\int_{\Rm_+}\int_{\Rm}\int_{\Rm} \int_{\Rm^3}\int_{\Rm^3}d q dk d|k'| d\theta^1_{k'}d\theta^1_{k'} d\xi d\xi'|k'|^2.
\end{align*}
Above, we used the notation $\xi_\eps=\xi_\pr \hat k+ \sqrt{\eps}\xi_\bot$ and $\xi'_\eps=\xi'_\pr \hat k+ \sqrt{\eps}\xi'_\bot$. The leading term is obtained for $\xi_\pr+\sigma_1 |\xi_\pr|=\xi'_\pr+\sigma_2 |\xi_\pr'|=0$ so that the first two phases vanish. There is otherwise some averaging that leads to a higher order contribution. Writing $q=q_\pr \hat k+ q_\bot$, $q_\bot \cdot k=0$,  integrating in $(\theta^1_{k'},\theta^2_{k'})$ the phase term involving $(e_1 \theta^1_{k'}+ e_2 \theta^2_{k'}) \cdot (\xi_\pr \hat k-2q)=-2(e_1 \theta^1_{k'}+ e_2 \theta^2_{k'}) \cdot q$, we obtain a Dirac measure enforcing that $q_\bot =0$. As a consequence, using the fact that $F(x)= F(|x|)$ and  $\hat \Phi(k)= \hat \Phi(|k|)$:
\begin{align*}
&V^C_{n,p} \sim \sigma_n^2 \eps^{4}e^{-\Sigma}\sum_{\pm} \int d Y F_\eps\big(||k|\mp q_\pr|\big) F_\eps\big(||k'|\pm q_\pr-\frac{1}{2} (|\xi_\pr|+|\xi'_\pr|)|\big) \\
&\hspace{3cm}\times e^{i |\xi_\bot|^2/(2|\xi_\pr|)}e^{i |\xi'_\bot|^2/(2|\xi'_\pr|)}e^{i z \cdot (\xi_\eps-\xi'_\eps)/\eps} \\
&\hspace{3cm}\times \hat \Phi\big(|\xi_\pr\mp 2q_\pr|\big) g_\mu(|\xi_\pr|-|k_0|)g_\mu (|\xi'_\pr|-|k_0|)\\
&\int dY=\int_{\Rm_+}\int_{\Rm_+}\int_{\Rm_+}\int_{\Rm^3}\int_{\Rm^3}\int_{\Rm^3}\int_{\Rm_+} dq_\pr d\xi_\pr d\xi_\pr' d\xi_\bot d\xi'_\bot dk d |k'| |k'|^2.
\end{align*}
The final expression is obtained by setting $|k| \to \eps^\gamma |k|/r_0\pm q_\pr$, 
$q_\pr\to \eps^\gamma q_\pr/r_0\mp |k'|\pm (|\xi_\pr|+|\xi'_\pr|)/2$, $|\xi_\pr| \to \mu^{-1}|\xi_\pr|+|k_0|$ and $|\xi'_\pr| \to \mu^{-1}|\xi'_\pr|+|k_0|$, which leads to, using that $\mu \gg 1$:
$$
V^C_{n,p}(z) \sim \sigma_n^2 \eps^{4(1-\gamma)} r_0^4 e^{-\Sigma} e^{i\frac{|z|^2}{\eps}}\int_{\Rm}\int_{\Rm} d\xi_\pr d\xi'_\pr J_0 \big(|z||(\xi_\pr-\xi_\pr')/(\mu \eps)\big)G(\xi_\pr,\xi'_\pr),
$$
for some smooth function $G$ and where $J_0$ is the zero-th order Bessel function of the first kind. Hence, $V^C_{n,p}(z) \sim \sigma_n^2 \eps^{4(1-\gamma)}r_0^4 e^{-\Sigma} $ for $|z| \leq \eps \mu$ and $V^C_{n,p}(z)\ll V^C_{n,p}(0)$ when $|z|\gg \eps \mu$ so that $V_{n,p}^C$ is essentially supported on $|z| \leq \eps \mu$. This ends the proof.

\subsection{Proof of proposition \ref{propvarK}} \label{proofpropvarK}
We use here the notation of section \ref{meas_imag}. With $\delta p^\eps$ the random fluctuations defined in \fref{randH}, the variance of the WB functional admits the expression
$$
V^W(z)=\frac{1}{(2 \pi)^3} \int_{\Rm^3} dk dk' e^{i (k-k') \cdot z} \cos |k|\cos |k'| \E \{ (\calF \un_D \delta p^\eps)(t=1,k)   \overline{(\calF \un_D \delta p^\eps)}(t=1,k')\}.
$$
Similarly to the CB functional, we neglect the finite size of the detector as a first approximation, so that $V^W$ reads
$$
V^W(z)\sim \int_{\Rm^3} dk dk' e^{i (k-k') \cdot z} \cos |k|\cos |k'| \E \{ (\calF \delta p^\eps)(t=1,k)    \overline{(\calF \delta p^\eps)}(t=1,k')\}.
$$
We start by computing $(\calF \delta p^\eps) (t,k)$. For a regular function $f$, we have
$$
(\Box^{-1} f)=\int_0^t G(t-s,\cdot) * f(s,\cdot) ds, \quad \mbox{so } (\calF (\Box^{-1} f)(t,k)=\int_0^t \frac{ds}{|k|} \sin |k| (t-s)\calF f(s,k) ds.
$$
Using \fref{balK} and \fref{randH}, this implies that
\bee
(\calF \delta p^\eps) (t,k)&=&-\sigma_0 e^{-t\Sigma/2}\int_0^t  \frac{ds}{|k|} \sin |k| (t-s) \calF \left(V \left(\frac{\cdot}{\eta \eps} \right)\frac{\partial^2 p^2_B}{\partial t^2}\right)(s,k)\\
&=&-\sigma_0 (\eta \eps)^3 e^{-t\Sigma/2}\int_0^t  \frac{ds}{|k|} \sin |k| (t-s) (\hat V(\eps \eta \cdot) *_k  \widehat{\Delta p^B})(s,k).
\eee
Recall that $p^B=\partial_t G(t,\cdot) * p_0^\eps$, so that $\Delta p^B=\partial_t G(t,\cdot) * \Delta p_0^\eps$ and
$$\widehat{\Delta p^B}(s,k)= -\eps^3|k|^2 g_\mu\left(\eps (|k|-\frac{|k_0|}{\eps})\right)\cos |k|s, \qquad g_\mu(x)=\mu g(\mu x).
$$
Hence
$$
(\calF \delta p^\eps) (t,k)=\sigma_0 \eta^3 \eps^6e^{-\frac{t\Sigma}{2}}\int_0^t \int_{\Rm^3} ds dp |k| \sin |k| (t-s)
\hat{V}(\eps \eta p) g_\mu\left(\eps (|k-p|-\frac{|k_0|}{\eps})\right) \cos |k-p|s.
$$
Therefore, using that $\E \{\hat{V}(\xi)\hat {V}(\nu) \}=  (2 \pi)^{3}\hat{R}(\xi) \delta_0(\xi+\nu)$, we find for the second moment
 \bal
 &\E \{(\calF \delta p^\eps) (t=1,k)\overline{(\calF \delta p^\eps)} (t=1,k')\}\sim \\
 &\sigma_0^2 (\eta \eps)^3 \eps^6\int_0^{1}\int_0^{1} \int_{\Rm^3} ds ds' dp|k||k'|   \sin |k| (1-s) \sin |k'| (1-s') f(k,k',p,s),
 \end{align*}
where
$$
f(k,k',p,s)=e^{-\Sigma}\hat{R}(\eps \eta p) g_\mu\left(\eps (|k-p|-\frac{|k_0|}{\eps})\right)  \overline{g_\mu}\left(\eps (|k'-p|-\frac{|k_0|}{\eps})\right) \cos |k-p|s \cos |k'-p|s'.
$$
Hence, the variance of the functional reads
\bea \nonumber
V^W(z)&\sim& \sigma_0^2 \eta^3 \eps^9 \int_0^1\int_0^1 \int_{\Rm^9} ds ds' dp dk dk'|k||k'|e^{i (k-k') \cdot z}\sin |k| (1-s) \sin |k'| (1-s') \\
&&\cos |k|\cos |k'|   f(k,k',p,s). \label{vartrd}
\eea
Setting in order $k' \to p+ \eps^{-1}k'$,  $k \to p+ \eps^{-1}k$, as well as $p \to \eps^{-1} p$, and writing the sines and cosines as sum of complex exponentials, we find

\bal
&V^W(z) \sim   \frac{ \sigma_0^2 \eta^3}{\eps^{2}} 
\sum_{\sigma_1, \sigma_2,\sigma_3,\sigma_4,\sigma_5,\sigma_6=\pm 1} \sigma_3 \sigma_4 \int dX |k+p| |k'+p| h_\mu(k',p,k)\\
&\hspace{1cm}  \times \exp \frac{i}{\eps}\left\{ \sigma_1 |k+p|+\sigma_2 |k'+p|\right\}\exp \frac{i}{\eps}\left\{ \sigma_3 | k+p|+\sigma_4 |k'+p|\right\}\\
& \hspace{1cm}\times \exp  \frac{i}{\eps}\left\{(\sigma_5 |k|-\sigma_3 |k+p|)s\right\}\exp \frac{i}{\eps}\left\{(\sigma_6 |k'|-\sigma_4 |k'+p|) s'\right\},
\end{align*}
where
$$
h_\mu(k',p,k)=e^{-\Sigma}e^{i (k-k')\cdot z/\eps}g_\mu(|k|-|k_0|)  \overline{g_\mu}(|k'|-|k_0|)\hat R(\eta p) 
$$
and
$$
\int dX=\int_0^{1}\int_0^1\int_{\Rm^3}\int_{\Rm^3} \int_{\Rm^3} ds ds' d k dk' dp. 
$$
The first two oscillating phases in $V^W$ compensate directly when $\sigma_1=-\sigma_3$ and $\sigma_2=-\sigma_4$. The phases are otherwise strictly positive which leads to some averaging. The leading order is therefore obtained for $\sigma_1=-\sigma_3$ and $\sigma_2=-\sigma_4$. We then use the following short time - long time decomposition of $V^W$:
\begin{align*}
&V^W(z) = \int_0^\eps\int_0^\eps ds ds'(\cdots)+ \int^1_\eps\int_\eps^1 ds ds'(\cdots)+\int_0^\eps\int^1_\eps ds ds'(\cdots)+\int_\eps^1\int_0^\eps ds ds'(\cdots)\\
&:=\sum_{i=1}^4 V_i(z).
\end{align*}
The terms $V_3$ and $V_4$ can be shown to be negligible compared to the first two because of the different integration domains. We focus consequently only on the most interesting terms $V_1$ since $V_2$. Let us start with the most technical term $V_2$. We perform a standard stationary phase analysis in the variable $p$. We are looking for a point $p_0$ such that $\sigma_3 \widehat{k+p_0} s+\sigma_4 \widehat{k'+p_0} s'=0$ (first order term in the phase), together with $\sigma_5 |k|-\sigma_3 |k+p_0|=\sigma_6 |k'|-\sigma_4 |k'+p_0|=0$ (zero order term). This suggests that $\sigma_5=\sigma_3$, $\sigma_6=\sigma_4$, $p_0=0$ and $\sigma_3 \widehat{k} s=-\sigma_4 \widehat{k'} s'$. Defining $\xi_1= s \hat k$, $\xi_2= s' \hat k'$, and performing the changes of variables $p \to \sqrt{\eps} p$, $\xi_1=-\sigma_3 \sigma_4 \xi_2+\sqrt{\eps} \xi_1$ leads to, all computations done for the leading term:
\bal
&V_2(z) \sim \sigma_0^2 \eta^{3-\delta} \eps^{1-\frac{\delta}{2}}e^{-\Sigma}
\sum_{\sigma_3,\sigma_4=\pm 1} \sigma_3 \sigma_4 \int dX^\eps |k|^3 |k'|^3 |\xi_1|^{-\delta} |\xi_2|^{\delta/2-4} e^{i\, \textrm{sign}(\sigma_3/(2|k|)+\sigma_4/(2|k'|)) |\xi_1|^2}  \\
& \hspace{7cm} \times |\sigma_3/(2|k|)+\sigma_4/(2|k'|)|^{\frac{\delta}{2}} |\sin (\hat \xi_1 \cdot \hat \xi_2)|^\delta\\
& \hspace{7cm}
\times e^{-i (\sigma_3 \sigma_4|k|+|k'|) \hat \xi_2 \cdot z/\eps}
g_\mu(|k|-|k_0|)  \overline{g_\mu}(|k'|-|k_0|)\\
&\int dX^\eps=\int_{\Rm^3} \int_{\eps \leq |\xi_2| \leq 1}\int_{\Rm_+} \int_{\Rm_+} d\xi_1 d\xi_2 d|k| d|k'|.
\end{align*}
Setting finally $\xi_2 \to \eps \xi_2$, $|k|\to |k_0|+ |k|/\mu$, $|k'|\to |k_0|+ |k'|/\mu$, we obtain
\bal
&V_2(z) \sim \sigma_0^2 \eta^{3-\delta} e^{-\Sigma}
 \int_{S^2} d\hat \xi_2 e^{-2i |k_0| \hat \xi_2 \cdot z/\eps} \left|\hat g \left(\frac{z \cdot \hat \xi_2}{\eps \mu}\right)\right|^2
\end{align*}
The function above has a very similar structure to one of the initial condition, and is therefore mostly supported on $|z| \leq \eps \mu$. 

Let us consider now the term $V_1$. After the changes of variables $s  \to \eps s$ and $s' \to \eps s'$, we find (recall that we have for the leading order that $\sigma_1=-\sigma_3$ and $\sigma_2=-\sigma_4$)
\bal
&V_1(z) \sim   \sigma_0^2 \eta^3 
\sum_{\sigma_3,\sigma_4,\sigma_5,\sigma_6=\pm 1} \sigma_3 \sigma_4 \int dX |k+p| |k'+p| h_\mu(k',p,k)\\
& \hspace{1cm}\times \exp  i\left\{(\sigma_5 |k|-\sigma_3 |k+p|)s\right\}\exp i \left\{(\sigma_6 |k'|-\sigma_4 |k'+p|) s'\right\}.
\end{align*}
When $\eta=1$, then $V_1$ is of the same order as $V_2$. When $\eta \ll 1$, we will see that $V_1$ is the leading term. Remarking first that
\bal
&\sigma_3 |k+p| \exp  i\left\{(\sigma_5 |k|-\sigma_3 |k+p|)s\right\}= \sigma_5 |k| \exp  i\left\{(\sigma_5 |k|-\sigma_3 |k+p|)s\right\}\\
& \hspace{7cm}- \frac{1}{i}\frac{d}{ds} \exp  i\left\{(\sigma_5 |k|-\sigma_3 |k+p|)s\right\},
\end{align*}
with a similar observation for the other exponential, $V^W_1$ can be split into four terms. It can be shown that the leading one is the one involving the product of the derivatives, so that
\bal
&V_1(z) \sim   \sigma_0^2 \eta^3 
\sum_{\sigma_3,\sigma_4,\sigma_5,\sigma_6=\pm 1} \int_{\Rm^9} dk dk' dp h_\mu(k',p,k)\\
& \hspace{1cm}\times \left(\exp  i\left\{(\sigma_5 |k|-\sigma_3 |k+p|)\right\}-1 \right)\left(\exp i \left\{(\sigma_6 |k'|-\sigma_4 |k'+p|)\right\}-1 \right).
\end{align*}
With the changes of variables $|k|\to |k_0|+ |k|/\mu$, $|k'|\to |k_0|+ |k'|/\mu$ and $p \to p/\eta$, we find
\bal
&V_1(z) \sim \sigma_0^2 e^{-\Sigma}
 \int_{S^2}\int_{S^2} d\hat k d \hat k' e^{i (\hat k-\hat k') \cdot z/\eps} \left|\hat g \left(\frac{(\hat k-\hat k') \cdot z}{\eps \mu}\right)\right|^2\\
&\hspace{3cm} \sum_{\sigma_3,\sigma_4=\pm 1}\int_{\Rm^3} dp \hat R(p) \left( e^{ - i \sigma_3 |p| /\eta)}-1 \right)\left(e^{-i \sigma_4 |p|)/\eta}-1 \right).
\end{align*}
Terms in the expression above involving an oscillating exponential are negligible, which shows that $V_1(0) \sim \sigma_0^2 e^{-\Sigma}$, and for the same reason as $V_2$, $V_1$ is mostly supported on $|z| \leq \eps \mu$. The leading term is therefore $V_1$, which ends the proof.
\subsection{Proof of proposition \ref{propvarKn}} \label{proofpropvarKn}
The proof is straightforward, the variance of the WB functional for the noise contribution admits the expression
$$
V^W_n(z)=\frac{\sigma_n^2}{(2 \pi)^3} \int_{\Rm^3} dk dk' e^{i (k-k') \cdot z} \cos |k|\cos |k'| \E \{ (\calF \un_D n_p^\eps)(t=1,k)    \overline{(\calF \un_D n_p^\eps)}(t=1,k')\}.
$$
Neglecting the finite size of the detector in first approximation leads to
$$
V^W_n(z)\sim \sigma_n^2 \eps^3 \int_{\Rm^3} dk   (\cos |k|)^2 \hat \Phi(\eps k) \sim \sigma_n^2.
$$
\section{Conclusion} \label{conc}
This work is concerned with the comparison in terms of resolution and stability of prototype wave-based and correlation-based imaging functionals. In the framework of 3D acoustic waves propagating in a random medium with possibly long-range correlations, we obtained optimal estimates of the variance and the SNR in terms of the main physical parameters of the problem. In the radiative transfer regime, we showed that for an identical cross-range resolution, the CB and WB functionals have a comparable SNR. The CB functional is shown to offer a better SNR provided the resolution is lowered, which is achieved by calculating correlations over a small domain compared to the detector. This is the classical stability/resolution trade-off. We obtained morever that the minimal central wavelength $\lambda_m$ that the functionals could accurately reconstruct were identical in the regime of weak fluctuations in the random medium, and that in the case of larger fluctuations, the CB functional offered a better (smaller) $\lambda_m$ (resolution) than the WB functional.

We also investigated the effects of long-range correlations in the complex medium. We showed that coherent imaging became difficult to implement because the mean free path was very small and therefore the measured signal too weak compared to additive, external, noise. In such a case, transport-based imaging with lower resolution is a good alternative to wave-based (coherent) imaging. This will be investigated in more detail in future works.
\section*{Acknowledgment}
This work was supported in part by grant AFOSR NSSEFF- FA9550-10-1-0194. 
{\footnotesize \bibliographystyle{siam}
  \bibliography{bibliography3} }

 \end{document}